\documentclass{amsart}
\usepackage{graphicx}
\title{A CLT for a band matrix model}
\author{Greg Anderson}
\address{School of Mathematics, University of Minnesota, 
206 Church St. SE,
Minneapolis, MN 55455}
\author{Ofer Zeitouni}
\thanks{ O.Z. was partially supported by NSF grant number
DMS-0302230.}
\address{School  of Mathematics, University of Minnesota, 
206 Church St. SE,
Minneapolis, MN 55455, and
Departments of Mathematics and of EE, Technion, Haifa 32000, Israel}
\date{June 15, 2004. Revised December 1, 2004.}

\newcommand{\RR}{{\mathbb{R}}}
\newcommand{\CC}{{\mathbb{C}}}
\DeclareMathOperator{\supp}{{\mathrm{supp}}}
\newcommand{\card}{{\#}}
\newtheorem{Theorem}[subsection]{Theorem}
\newtheorem{Assumption}[subsection]{Assumption}
\newtheorem{Lemma}[subsection]{Lemma}
\newtheorem{Definition}[subsection]{Definition}
\newtheorem{Proposition}[subsection]{Proposition}

\DeclareMathOperator{\trace}{{\mathrm{tr}}}
\DeclareMathOperator{\weight}{{\mathrm{wt}}}
\DeclareMathOperator{\Var}{{\mathrm{Var}}}
\newcommand{\Poincare}{{Poincar\'{e} }}
\newcommand{\FFF}{{\mathcal{F}}}
\newcommand{\NNN}{{\mathcal{N}}}
\newcommand{\one}{{\mathbf{1}}}

\DeclareMathOperator{\Res}{{\mathrm{Res}}}
\begin{document}
\bibliographystyle{alpha}
\maketitle
\begin{abstract}
A law of large numbers and a
central limit theorem are derived for linear statistics of
random symmetric  matrices whose on-or-above diagonal entries
are independent, but neither necessarily identically distributed,
nor necessarily all of the same variance. The
derivation is based on systematic combinatorial enumeration, study of
generating  functions, and concentration inequalities of the
\Poincare type.  Special cases  treated,  with an explicit evaluation
of limiting variances,  are generalized Wigner
and Wishart matrices.
\end{abstract}
\section{Introduction}
\label{section:Intro}
The interest in the limiting properties of the empirical distribution
of eigenvalues of large symmetric random  matrices can be traced back to
\cite{wishart} and to the
path-breaking article of Wigner \cite{wigner}.
We refer to \cite{bai}, \cite{deift}, \cite{petz}, \cite{mehta} 
and \cite{pastur}
for partial overview and some of the recent spectacular progress
in this field.

In this paper we study both convergence of the empirical 
distribution
and central limit theorems for linear statistics of the empirical
distribution of a class of random matrices. To give right away the flavor
of our results, 
consider for each positive integer $N$ the $N$-by-$N$ symmetric random
matrix
$X(N)$ with on-or-above-diagonal entries
$X(N)_{ij}=N^{-1/2} f(i/N,j/N)^{1/2} \xi_{ij}$, where
the $\xi_{ij}$ are zero mean unit variance i.i.d.\ random variables
satisfying the \Poincare inequality with constant $c$ (see 
\S\ref{subsec-poincarematrices} for 
the definition), and $f(\cdot,\cdot)$ is a  nonnegative function 
symmetric and
continuous on $[0,1]^2$ such that 
$\int_0^1 f(x,y)dy\equiv 1$. 
Define the {\it semicircle distribution} $\sigma_S$ of zero mean and
unit variance  to be the measure
on $\RR$ of compact support with density
$\frac{d\sigma_S}{dx}:=\frac{1}{2\pi} \sqrt{4-x^2}{\bf 1}_{\{|x|\leq
2\}}$\,. Let $\lambda_1(N)\leq \cdots\leq \lambda_N(N)$ be the
eigenvalues of $X(N)$.
Under  these assumptions, 
a corollary of our general results (see 
Theorem \ref{Theorem:Mainspecial} below) states that the 
{\it empirical distribution} $L(N):=N^{-1} \sum_{i=1}^N
\delta_{\lambda_i(N)}$ converges weakly, in probability, to $\sigma_S$,
and further, for any  continuously differentiable function $f$ on $\RR$ of
polynomial growth, with $\|f'\|_{L^2(\sigma_S)}>0$,
the sequence of random variables
\begin{equation}
\label{eq:040204a}
 \sum_{i=1}^N f(\lambda_i(N))-E \left[ \sum_{i=1}^N
f(\lambda_i(N))\right]
\end{equation}
converges in distribution to a nondegenerate zero mean Gaussian random
variable with variance given by an explicit formula. Similar explicit
results hold for a class of generalized Wishart matrices, see Theorem
\ref{theo-groupwishart}.  For polynomial test functions
$f$, such results hold in much greater generality, see Theorem 
\ref{Theorem:Main}.

Our approach has two main components. The first, of some interest on its
own,
is a combinatorial
enumeration scheme for the different types of terms that  contribute to
the expectation of products of traces of powers of the matrices under
study.
This scheme takes the bulk of the paper to develop.
The other component, which allows us to move from polynomial test
functions to continuously differentiable ones, is based on concentration
inequalities of the \Poincare type. The latter component is developed in
\S\ref{sec-poincareconc}, building on earlier results
of concentration for random matrices that can be found in
\cite{Chabo} and
\cite{GZconc}.

CLT results related to 
Theorems \ref{Theorem:Main}
and \ref{Theorem:Mainspecial} below
have already been stated in the literature. An especially strong
inspiration
to our study is the work of Jonsson \cite{jonsson},
who gives CLT statements for traces of polynomial functions of Gaussian
Wishart matrices, based on the method of moments
introduced by Wigner in \cite{wigner}.  The method of moments was
revisited in the far-reaching work of Sinai and Soshnikov
\cite{soshnikov}, where, as
an easy by-product of their results, they state a CLT for 
traces of analytic functions of
Wigner-type matrices.
Pastur and co-authors, on the one hand,
and Bai and co-authors, on the other, 
have championed an approach based on the
evaluation of resolvents. The latter approach has the  advantage of
allowing one to relax hypotheses on matrix entries; in  particular one does
not need to have all moments finite. CLT statements based on these
techniques and expressions for the resulting variance, for functions of the
form
$f(x)=\sum a_i/(z_i-x)$ where $z_i\in\CC\setminus \RR$, and matrices 
of Wigner type, 
can be found in \cite{khorunzhy1}, with somewhat sketchy proofs. Earlier 
statements can be found in \cite{girko}. A complete
treatment for $f$ analytic in a domain including the support
of the limit 
of the empirical distribution of eigenvalues is given in \cite{bainew} 
for matrices  of Wigner type, and in 
\cite{bai-silvernew} for matrices
of Wishart type
under a certain restriction on fourth moments. Much more is known for
restricted classes of matrices:
Johansson
\cite{johansson},
using an approach based
on the explicit joint density of the eigenvalues available in the
independent case 
only in the Gaussian Wigner situation, 
characterizes
completely those functions $f$ for which a CLT holds. Cabanal-Duvillard
\cite{CD01}
introduces a stochastic calculus approach and 
proves a CLT for traces of polynomials of
Gaussian Wigner and Wishart matrices, as well as for traces of
non-commutative 
polynomials of pairs of independent Gaussian Wigner matrices. 
Recent extensions and reinterpretation of his work,
using the notion of second order freeness, can be found in
\cite{mingospeicher}. Still in the
Gaussian case, Guionnet \cite{guionnet-clt}, using a stochastic calculus
approach, gives  a CLT (with a somewhat implicit variance computation)
for  a class of functions $f$ in the case of {\it band matrices}. 
Earlier, laws of large
numbers for band matrices were derived, see e.g.
\cite{MPK}, \cite{Sh96} and the references therein.
In comparison with
the references mentioned above, our work can be seen as relaxing the 
structural assumptions
on the variance of the entries of the matrix $X(N)$, as well as the
Gaussian assumption, while 
still requiring rather strong moment bounds on
the individual entries (if one is interested only in polynomial test
functions) or \Poincare type conditions on the entries (if one wants a
wider class of test functions). 

The structure of the article is as follows. In \S\ref{section:Model},
we introduce the matrix model considered throughout the paper and set
basic notations. Our main results for polynomial test functions $f$ are
stated 
in \S\ref{Section:maintheorem}. \S\ref{section:grammar}
develops the language we use in the combinatorial enumeration mentioned
above.
\S\ref{section:Limcal} is devoted to some preliminary
limit calculations which are then immediately applied in \S\ref{sec:proofI}
to prove our main result 
concerning limiting spectral measures, Theorem~\ref{Theorem:Preliminary}.
\S\ref{sec:FK} is devoted to the derivation of some {\em a priori}
estimates, following \cite{furedi}, useful in the study of the support of
the empirical distribution 
$L(\NNN)$. \S\ref{sec:bracelets} is  the heart of our enumeration scheme,
and the results are  immediately applied in \S\ref{sec:proofII} to yield
the proof of our main CLT statement, Theorem~\ref{Theorem:Main}.
\S\ref{sec:CalcMean} is devoted to the proof of
Theorem~\ref{Theorem:MainBis}, which is a technical result
describing how to approximate
$E\trace X(\NNN)^n$ at CLT scale; this part of the paper may be skipped
without much loss of comprehension of the remainder of the paper.
\S\ref{sec-poincareconc} is devoted to concentration of measure results
based on the \Poincare inequality. Finally, in
\S\ref{sec-ChebDiag} we specialize our main results to generalized Wigner
and Wishart matrices, and derive explicit representations for the
resulting variances.

\section{The model}
\label{section:Model}
We define in this section the class of random matrices we are going to
deal with. Matrices of this class are symmetric, with
on-or-above-diagonal entries independent, with all entries possessing
moments of all orders, and with
off-diagonal entries of mean zero; further and crucially, subject to the
constraints of symmetry and vanishing of off-diagonal means, the moments of
entries of such matrices are allowed to depend upon position. Now as it
turns out, only certain statistical properties of the patterns of first,
second and fourth moments of entries figure in our limit formulas. 
Accordingly, our description of the class is contrived so as to  emphasize
those statistical properties and to suppress unneeded detail concerning the
exact dependence of moments of entries on position.  The notion crucial for
gaining ``statistical control'' is that of {\em color}. The reader
interested only in Wigner matrices should take as space of colors a space
consisting of a single color. 

\subsection{The band matrix model}\label{subsection:BandStructure}
\subsubsection{Colors}\label{subsubsection:ColorHypotheses}
We fix a Polish space, elements of which we call {\em colors}.
We declare the Borel sets of color space to be measurable.   We fix a
probability  measure $\theta$ on color space. We fix a bounded measurable 
real-valued function
$D$ on color space.  For each positive integer $k$ we fix 
a bounded measurable nonnegative function $d^{(k)}$ on color
space and a symmetric bounded measurable nonnegative function
$s^{(k)}$ on the product of two copies of color space.
We make the following assumptions:
\begin{itemize}
\item $d^{(k)}$ is constant for $k\neq 2$.
\item $s^{(k)}$ is constant for $k\not\in\{2,4\}$.
\item $s^{(k)}$ has discontinuity set of measure zero with respect to
$\theta\otimes \theta$.
\item 
$D$, $d^{(k)}$ and the diagonal restriction of $s^{(k)}$ have
discontinuity sets of measure zero with respect to $\theta$. 
\end{itemize}
For any
bounded function
$f$ on color space, or on a product of copies of color space, we write
$|f|_\infty$ for its supremum norm.  

\subsubsection{Letters}
We fix a countably  infinite set, elements of which
we call {\em letters}.  
We fix a
function
$\kappa_0$ from letter space to color space, and we say that
$\kappa_0(\alpha)$ is the {\em color} of  the letter $\alpha$.
 Given any nonempty finite  set $\NNN$ of letters of
cardinality $N$ put
$$\theta_\NNN= N^{-1}\sum_{\alpha\in \NNN}\delta_{\kappa_0(\alpha)},$$
which  is the color distribution of letters belonging to $\NNN$.
We reserve the script letter $\NNN$ for use in this context and invariably
denote the cardinality of $\NNN$ by the roman letter $N$. Analogously,
given a sequence
$\NNN_1,\NNN_2,\NNN_3,\dots$ of finite nonempty sets of letters,
$N_1,N_2,N_3,\dots$ denotes the corresponding sequence of cardinalities.

\subsubsection{The family $\{\xi_e\}$ of random variables}
We fix a family 
$\{\xi_e\}$ of independent real-valued mean zero random variables
indexed by unordered pairs $e$ of letters. We
assume that for all letters
$\alpha,\beta$ and positive integers $k$  we
have
\begin{equation}\label{equation:MomentConstraint}
E|\xi_{\{\alpha,\beta\}}|^k\leq
\left\{\begin{array}{cl}
s^{(k)}(\kappa_0(\alpha),\kappa_0(\beta))&\mbox{if $\alpha\neq \beta$,}\\
d^{(k)}(\kappa_0(\alpha))&\mbox{if $\alpha=\beta$,}
\end{array}\right.
\end{equation}
and moreover we assume that  equality holds above whenever
one of the following conditions holds:
\begin{itemize}
\item $k=2$.
\item
$\alpha\neq \beta$ and $k=4$. 
\end{itemize}
In other words, the rule is to enforce
equality whenever the not-necessarily-constant functions
$d^{(2)}$,
$s^{(2)}$ or
$s^{(4)}$ are involved, but otherwise merely to impose a bound.

\subsubsection{Random matrices}
 Given any nonempty finite set
$\NNN$ of letters, let
$X(\NNN)$ be the $N\times N$  real symmetric random matrix with
entries 
$$X(\NNN)_{\alpha\beta}
=D(\kappa_0(\alpha))\delta_{\alpha\beta}+
N^{-1/2}\xi_{\{\alpha,\beta\}}\;\;\;\;\;\;(\alpha,\beta\in \NNN),$$ 
denote the eigenvalues of $X(\NNN)$ by
$\lambda_1(\NNN)\leq\cdots\leq\lambda_{N}(\NNN)$,
 and let
$$L(\NNN)= N^{-1}\sum_{i=1}^{ N} \delta_{\lambda_i(\NNN)}$$
be the empirical
distribution of the spectrum 
of $X(\NNN)$. Put
$$\overline{L}(\NNN)=EL(\NNN).$$
Note that
$$\langle \overline{L}(\NNN),x^n\rangle=\frac{1}{N}E\trace X(\NNN)^n,$$
where here and often below we employ the abbreviated notation 
$$\langle \mu,f\rangle=\int
f(x)\mu(dx)\,$$ 
for integrals.

\subsection{Generating functions}
\subsubsection{}
Let $\sigma$ be any probability measure on color space.
Let $$[\Phi_{n,\sigma}(c)]_{n=1}^\infty$$ be the unique sequence
of real-valued bounded measurable functions on color space characterized
by the generating function identity
\begin{equation}\label{equation:PhiRecursion}
\Phi_\sigma(c,t)=\left(\frac{t}{1-D(c)t}\right)
\left(1-\frac{t}{1-D(c)t}\int
s^{(2)}(c,c')\Phi_\sigma(c',t)\sigma(dc')\right)^{-1}
\end{equation}
where
$$\Phi_\sigma(c,t)=\sum_{n=1}^\infty \Phi_{n,\sigma}(c)t^n$$
is the corresponding generating function. We emphasize that we
view the power series here formally, i.~e., as devices for managing
sequences, not as analytic functions. We write
(\ref{equation:PhiRecursion}) as a shorthand for the
recursion obtained by formally expanding  both sides of
(\ref{equation:PhiRecursion}) in powers of $t$,  and then equating
coefficients of like powers of $t$. When $\sigma=\theta$, we omit 
it from the notation.

\subsubsection{}
For each positive integer $r$ we define a function
$$
K_r(c_1,\dots,c_r)=\left\{\begin{array}{cl}
s^{(2)}(c_1,c_1)&\mbox{if $r=1$,}\\
s^{(2)}(c_1,c_2)^2&\mbox{if $r=2$,}\\
s^{(2)}(c_1,c_2)s^{(2)}(c_2,c_3)\cdots
s^{(2)}(c_r,c_1)&\mbox{if
$r\geq 3$,}
\end{array}\right.
$$
on the product of $r$ copies of color space. We define
\begin{equation}\label{equation:ThetaDef}
\Theta(x,y)=\sum_{r=1}^\infty
\frac{1}{r}\int\cdots \int K_r(c_1,\cdots,c_r)
\prod_{i=1}^r\left(\Phi(c_i,x)\Phi(c_i,y)\theta(dc_i)\right).
\end{equation}
We view $\Theta(x,y)$ as a formal power series in $x$
and $y$ with real coefficients; in keeping with this point of view, the
integrals on the right side of (\ref{equation:ThetaDef}) are  to be
evaluated by  first expanding the integrands in powers of $x$ and
$y$ and then integrating term by term (and hence, all integrals, being
expectations of bounded measurable functions, are well defined). 

\subsubsection{}
Put
\begin{equation}\label{equation:Correction}
\begin{array}{rcl}
\Psi(x,y)&=&\displaystyle\int
(d^{(2)}(c)-2s^{(2)}(c,c))\Phi(c,x)\Phi(c,y)\theta(dc)\\\\
&&\displaystyle +\frac{1}{2}\int\int
(s^{(4)}(c_1,c_2)-3s^{(2)}(c_1,c_2)^2)\\\\
&&\displaystyle\;\;\;\;\;\;\times\Phi(c_1,x)\Phi(c_2,x)
\Phi(c_1,y)\Phi(c_2,y)
\theta(dc_1)\theta(dc_2).
\end{array}
\end{equation}
We view $\Psi(x,y)$ as a formal power series in $x$ and $y$ with
real coefficients. As above, the integrals 
are to be evaluated by first expanding integrands in powers of
$x$ and $y$ and then integrating term by term.

\subsubsection{}
In order to gain convenient access to the information coded in the
formal power series
$\Theta(x,y)$ and $\Psi(x,y)$ we introduce
the following (abuse of) notation. We write
$$\left\langle \sum_{i=0}^\infty a_it^i,
\sum_{i=0}^\infty b_jt^i\right\rangle=\sum_{i=0}^\infty
a_ib_i$$ for any sequences $[a_i]_{i=0}^\infty$
and $[b_j]_{j=0}^\infty$ of real numbers such that the sum on the right
has only finitely many nonzero terms. Similarly we write
$$\left\langle \sum_{i=0}^\infty \sum_{j=0}^\infty
a_{ij}x^iy^j,
\sum_{i=0}^\infty \sum_{j=0}^\infty b_{ij}x^iy^j\right\rangle
=\sum_{i=0}^\infty \sum_{j=0}^\infty a_{ij}b_{ij}$$
for any doubly infinite sequences
$[a_{ij}]_{i,j=0}^\infty$ and $[b_{ij}]_{i,j=0}^\infty$
of real numbers such that the sum on the right has only finitely many
nonzero terms. 

The following fact, proved in \S\ref{sec:proofI},
 explains the role of the sequence
$\Phi_{n,\sigma}$:
\begin{Lemma}\label{LemmaTheorem:Preliminary}
If $\sigma=\theta$ or $\sigma=\theta_\NNN$ for some finite nonempty set
of letters $\NNN$, then there exists a unique
probability measure
$\mu_\sigma$ on the real line such that
\begin{equation}\label{equation:PhiMomentFormula}
\langle \mu_\sigma, x^n\rangle=\langle
\sigma,\Phi_{n+1,\sigma}\rangle\;\;\;(n=0,1,2,\dots),
\end{equation}
\begin{equation}\label{equation:MuSupportBound}
\supp \mu_\sigma\subset [-C,C]\;\;\;\;(C=2(|D|_\infty+
|s^{(2)}|_\infty^{1/2})).
\end{equation}
\end{Lemma}
\noindent 
In what follows, we write $\mu=\mu_\theta$
and
$\mu_\NNN=\mu_{\theta_\NNN}$.
\section{Assumptions and main theorem}
\label{Section:maintheorem}
Throughout, we let
$\Rightarrow$ denote the weak convergence of probability measures.
All of our results will be obtained under the following basic
\begin{Assumption}\label{ass-main} 
All the assumptions and notations in \S\ref{section:Model} hold. 
Further, there exists a sequence
$[\NNN_k]_{k=1}^\infty$ of finite nonempty sets of letters such that $
N_k\rightarrow\infty$ and $\theta_{\NNN_k}\Rightarrow \theta$.
\end{Assumption}
\noindent Our main results are:
\begin{Theorem}\label{Theorem:Preliminary}
Let Assumption \ref{ass-main} hold. Then:
 (i)
$L(\NNN_k)\Rightarrow
\mu$ in probability. \linebreak (ii) $\mu_{\NNN_k}\Rightarrow\mu$.
\end{Theorem}
\begin{Theorem}\label{Theorem:Main}
Let Assumption \ref{ass-main} hold. Fix
a real-valued polynomial function $f(\cdot)$ on the real line. Then
the sequence of random variables
$$Z_{f,k}:=\trace f(X(\NNN_k))-E\trace f(X(\NNN_k))$$
converges in distribution  to a zero mean Gaussian 
random variable
$Z_f$ of variance
\begin{equation}
\label{equation:CovarianceFormula}
EZ_f^2=\langle 2\Theta(x,y)+\Psi(x,y),xf'(x)yf'(y)\rangle.
\end{equation}
\end{Theorem}
\begin{Theorem}\label{Theorem:MainBis}
In the setting of the preceding theorem, we also have
\begin{equation}\label{equation:MeanFormula}
\lim_{k\rightarrow\infty} E\trace
f(X(\NNN_k))-
N_k\cdot\langle\mu_{\NNN_k},f\rangle=\frac{1}{2}\langle
\Theta(t,t)+\Psi(t,t),tf'(t)\rangle=:E_f.
\end{equation}
\end{Theorem}
\noindent 
We state the formulas
(\ref{equation:CovarianceFormula}) and (\ref{equation:MeanFormula})
in separate theorems because their proofs are separated in the main body
of the paper. In fact, we have structured the paper so that the reader
interested only in (\ref{equation:CovarianceFormula}) and its
applications can largely ignore the extra  (and somewhat heavy)
apparatus needed to prove (\ref{equation:MeanFormula}).

The results above can be made more 
transparent, and their range extended,
for certain special cases.  Of particular interest is the following:
\begin{Theorem}\label{Theorem:Mainspecial}
Let Assumption \ref{ass-main} hold, 
and further assume that
\begin{equation}
\label{eq-271004}
D\equiv 0,\;\;\;\int s^{(2)}(c,c')\theta(dc')\equiv 1.
\end{equation}
Then: (i) $\mu$ is the semicircle law $\sigma_S$ 
of zero mean and unit variance. (ii)
For polynomial functions $f$ the random variables
$Z_{f,k}$ converge in distribution toward a mean zero Gaussian random
variable
$Z_f$ with variance given by 
(\ref{equation:CovarianceFormulaWignerRewrite})\,.
(iii) If the random variables $\xi_{\{\alpha,\beta\}}$ all
satisfy a
\Poincare  inequality with common constant $c$ (see
\S\ref{sec-poincareconc}  for definitions),
then statement (ii) extends to continuously differentiable 
functions $f$  with polynomial growth, with variance again given by 
(\ref{equation:CovarianceFormulaWignerRewrite})\,.
\end{Theorem}
\noindent We refer to  the situation 
in Theorem \ref{Theorem:Mainspecial}
above as the
{\em generalized
Wigner matrix} model, 
because when $s^{(2)}\equiv 1$ one recovers Wigner matrices.
The expression $E_f$ in  (\ref{equation:MeanFormula}) 
can also be computed in
this case, see (\ref{eq:300304b}) below. 
We note in passing that for Gaussian matrices,
the condition (\ref{eq-271004}) has  been identified 
in \cite[Corollary 3.4]{NSR02} as sufficient and necessary (if $D=0$)
for $\mu$ to equal the semicircle distribution.

Similar considerations apply to  the
{\em generalized
Wishart matrix} model, see 
\S\ref{subsec-secondspec} for details.

\section{Basic spelling, grammar and counting}
\label{section:grammar}
We introduce in this section the basic language employed throughout
the paper for discussing enumeration problems. From
letters we build {\em words}, from words we build {\em sentences}, and then
we distinguish certain classes of words and sentences in terms of
properties of naturally associated {\em graphs}.  The classes of words and
sentences singled out here for special attention are eventually 
going to be used to enumerate the
 terms in sums giving the (mixed and/or centered) moments of traces of
powers of our random matrices. In particular, the {\em Wigner words}
enumerate the only terms whose contributions to the law of large numbers
for linear statistics do not vanish in the limit, whereas the {\em CLT
word-pairs} take care of the only terms whose contributions to the CLT
variance do not vanish in the limit. Further, the {\em CLT sentences}
(which can be built up systematically from the CLT word-pairs) enumerate
the nonnegligible terms in sums giving mixed centered 
moments of traces of
powers of our random matrices. Critical weak Wigner words, and 
marked Wigner words, are needed (only) in the evaluation of the mean
shift of linear statistics.

\subsection{Words and sentences}
A {\em word} is a finite sequence of letters at least one letter long.
(Words are never empty!)
We denote the length of a word $w$ by $\ell(w)$.  We say that a word $w$ is
{\em closed} if the first and last letters of $w$ are the same. (Every
one-letter word is automatically closed.) We view letters as one-letter
words.
A {\em sentence} is a finite sequence of words at least one word long.
(Sentences are never empty, nor do they contain empty words!)
We view words as one-word sentences.
The {\em support}
$\supp a$ of a sentence $a$ is the set of letters appearing in $a$,
and the {\em combinatorial weight} $\weight a$ is
the cardinality of $\supp a$. 
 We say
that sentences
$a$ and
$b$ are {\em disjoint} if $\supp a\cap \supp b=\emptyset$.
We say that sentences $a$ and $b$ are {\em
equivalent} and write
$a\sim b$ if there exists a one-to-one letter-valued function $\psi$
defined on $\supp a$ such that the result of applying
$\psi$ letter by letter to
$a$ is $b$. In other words, $a\sim b$ whenever
 $a$
codes to
$b$ under a simple substitution cipher.

We warn the reader that we distinguish between a sentence $a$ and the word
$w$ obtained
by concatenating all words in $a$; the ``punctuation'' carries
information important for our purposes and therefore must not be ignored.
For example, taking the set $\{1,2,3\}$ temporarily as our alphabet,
the word $123123$, the two-word sentence $[123,123]$
and the three-word sentence $[1,231,23]$ are distinct objects according
to our point of view. 

\subsection{Graphs} 
We fix terminology concerning graphs in a slightly
restrictive but convenient way as follows.  A {\em graph}
$G=(V,E)$ is an ordered pair consisting of a finite nonempty set $V$ of
letters and a set $E$ (possibly empty),  where each element of $E$ is
an unordered pair of elements of
$V$, i.~e., a subset of $V$ of cardinality $1$ or $2$. Elements of
$V$ are called {\em vertices} of
$G$, elements of
$E$ are called {\em edges} of $G$, and edges of cardinality $1$ are said
to be {\em degenerate}.  We say that a word
$w=\alpha_1\cdots \alpha_n$ of $n$ letters is a {\em walk} on $G$ provided
that
$\alpha_i\in V$ for $i=1,\dots,n$, and $\{\alpha_{i},\alpha_{i+1}\}\in E$
for
$i=1,\dots,n-1$, in which case we say that each of the vertices
$\alpha_i$ and edges
$\{\alpha_{i},\alpha_{i+1}\}$ of $G$ is {\em visited} by $w$.  A {\em
geodesic} in
$G$ is a walk visiting no vertex more than once. We say that
$G$ is {\em connected} if any two vertices are joined by a walk.
If $G$ is connected then $\card E\geq \card V-1$.  We call
$G$ a {\em tree} if
$G$ is connected  and $G$ has no nontrivial loops (in particular,
$G$ has no degenerate edges). Every two vertices of a tree are joined by a
unique
geodesic. For $G$ to be a tree it is necessary and sufficient that $G$ be
connected and 
$\card E\leq \card V-1$.
 A graph
$G'=(V',E')$ where
$V'\subset V$ and
$E'\subset E$ is called a {\em subgraph} of $G$.  A {\em
connected component} of $G$ is a connected subgraph of $G$ maximal in the
family of connected subgraphs of $G$. We call $G$ a {\em
forest} if every connected component of $G$ is a tree.
A {\em spanning forest} in $G$ is a graph $G'=(V',E')$
with $V'=V$ and $E'\subset E$ such that $G'$ is a forest having the
same number of connected components as does $G$. Every graph contains
at least one spanning forest.

\subsection{Orthographic and grammatical notions}
\subsubsection{The graph associated to a sentence}
Given a sentence 
$$a=[w_i]_{i=1}^n=[[\alpha_{ij}]_{j=1}^{\ell(w_i)}]_{i=1}^n$$ 
consisting of $n$ words (following a pattern we use often in
the sequel, $w_i$ denotes the
$i^{th}$ word of the sentence, and $\alpha_{ij}$ denotes the
$j^{th}$ letter of the $i^{th}$ word) we define 
$$G_a=(V_a,E_a)$$ to be the graph with
$$V_a=\supp a,\;\;\;E_a=\left\{\{\alpha_{ij},\alpha_{i,j+1}\}
\left|\begin{array}{l}
i=1,\dots,n,\\
j=1,\dots,\ell(a_i)-1
\end{array}\right.\right\}.$$
We view each word $w_i$ of the sentence $a$ in the natural way as a walk
on
$G_a$. 
We emphasize that
$E_a=\emptyset$ if $a$ consists of one-letter
words.
Note also  the difference between the graph associated to the sentence
$a$ and the graph associated to the single word consisting of the
concatenation of the words of $a$; in general the former has fewer edges
than the latter.

\subsubsection{Weak Wigner words} 
A word $w$ is called a {\em weak Wigner word} 
under the following two conditions:
\begin{itemize}
\item 
$w$ is closed.
\item $w$ 
visits every edge of $G_w$ at least twice.
\end{itemize}
Suppose now that $w$ is a weak Wigner
word. If
$\weight w={(\ell(w)+1)}/{2}$, then we drop the modifier ``weak'' and
call
$w$ a {\em Wigner word}.  (Every single letter word is automatically a 
Wigner word.)
If $\weight
w={(\ell(w)-1)}/{2}$, then we call
$w$  a {\em critical} weak Wigner word. 
For example, spelling with the alphabet $\{1,2,3\}$, 
we have that $w=121$  is a Wigner word and that
$w=12121$ is a critical weak Wigner word.

\subsubsection{Weak CLT sentences}
Let $a=[w_i]_{i=1}^n$
be a sentence consisting of $n$ words.
We say that $a$ is a {\em weak CLT sentence}
under the following three conditions: 
\begin{itemize}
\item All the words $w_i$ are closed.
\item Jointly the words/walks $w_i$ visit each edge of $G_a$ at least two
times.
\item For each $i\in \{1,\dots,n\}$ there exists $j\in
\{1,\dots,n\}\setminus\{i\}$ such that the graphs
$G_{w_i}$ and $G_{w_j}$ (both of which are subgraphs of $G_a$) have an
edge in common.
\end{itemize}
 Suppose now that
$a$ is a weak CLT sentence.  If
$\weight a=\sum_{i=1}^n \frac{\ell(w_i)-1}{2}$,
then we drop the modifier
``weak'' and call $a$ a {\em CLT sentence}.
If $n=2$ and $a$ is a CLT sentence, then we call $a$ a {\em CLT
word-pair}. 
For example, again spelling with the alphabet $\{1,2,3\}$,
we have that  $a=[1231,1321]$ is a CLT word-pair.

\subsubsection{Marked Wigner words}
A {\em marked Wigner word} is a three-word sentence $[w,\alpha,\beta]$
where $w$ is a Wigner word, and $\alpha$ and $\beta$ are distinct letters
appearing in $w$.

\subsubsection{Cyclic
permutations}\label{subsubsection:CyclicPermutations} Given a word
$w=[\alpha_i]_{i=1}^n$ of length $n$ and a permutation
$\sigma$ of $\{1,\dots,n\}$, we define $w^\sigma$ to be the word
$[\alpha_{\sigma(i)}]_{i=1}^n$. If $\sigma$ is a power of the
cycle $(123\cdots n)$, then we call
$\sigma$ a {\em cyclic permutation}, and we say that $w^\sigma$ is a {\em
cyclic permutation} of $w$.

\begin{Lemma}[``The parity principle''] Let $G$ be a forest. Let
$e$ be an edge of
$G$. Let $w$ be a word admitting interpretation as a  walk on $G$. Let
$w_*$ be the unique geodesic in
$G$ with initial and terminal vertices coinciding with those of
$w$. 
Then the word/walk
$w$ visits the edge
$e$ an odd number of times if and only if the geodesic $w_*$ visits $e$.
\end{Lemma}
\noindent This simple principle is repeatedly applied in the sequel. The
proof is elementary and therefore omitted.

\begin{Proposition}\label{Proposition:ExtremalWeakWigner}
Let $w$ be a weak Wigner word.
(i) We have $\weight
w\leq\frac{\ell(w)+1}{2}$ with equality
if and only if
$G_w$ is a tree. (ii) If $\weight w=\frac{\ell(w)+1}{2}$
then $w$ visits every edge of the tree $G_w$ exactly twice. (iii) $w$ is a
Wigner word if and only if there exists a decomposition
$w=\alpha w_1\cdots
\alpha w_r\alpha$ where $\alpha$ is the first letter of
$w$
and $w_1,\dots,w_r$ are pairwise disjoint Wigner words in which
$\alpha$ does not occur. 
 (iv)  The inequality $\frac{\ell(w)-1}{2}<\weight
w<\frac{\ell(w)+1}{2}$ is impossible. 
\end{Proposition}
\noindent These are ideas coming up in some proofs of Wigner's semicircle
law by  the method of moments. 
\proof Put $G=(V,E)=G_w=(V_w,E_w)$.
 (i)  The existence of the walk $w$ makes it clear that $G$ is
connected.  We have
\begin{equation}\label{equation:IntermediateInequality}
\weight w-1\leq \card E\leq \frac{\ell(w)-1}{2},
\end{equation} 
on the left because $G$ is connected, and on the right
by the hypothesis that $w$ is a weak Wigner word. The result follows.
(ii) Clear. 
(iii)($\Leftarrow$) Trivial.
(iii)($\Rightarrow$) By (i) and (ii) already proved,
the parts
of the tree
$G$ explored by the walk
$w$ between successive visits to the vertex $\alpha$ have to be disjoint.
(iv) Suppose rather that the inequality in question holds.
Then
$\ell(w)$ is even and
$\card V=\frac{\ell(w)}{2}$, hence by
(\ref{equation:IntermediateInequality}) we have
$\card E=\frac{\ell(w)}{2}-1=\card V-1$, and hence
$G$ is a tree. We now arrive at a contradiction:
by the parity principle $w$ cannot be both closed and a walk
that takes an odd number of steps.
\qed

As a consequence of Proposition \ref{Proposition:ExtremalWeakWigner},
one can visualize equivalence classes of Wigner words as rooted planar
trees, 
with the Wigner word determining an exploration  path on the tree that
visits
each vertex at least once and goes over each edge exactly twice, c.f.\
Figure 1. We do not make explicit use of this correspondence but it does
drive much of our intuition.
\begin{figure}[h]
\begin{center}
\includegraphics[scale=.49]{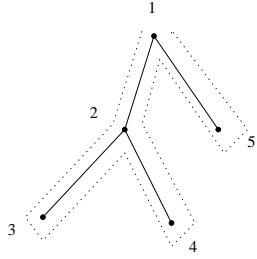}
\caption{The rooted planar tree (solid) and exploration process (dotted)
corresponding to the equivalence class of the Wigner word
$w=123242151$. Note the decomposition $w=1w_11w_21$ with $w_1=23242$ and
$w_2=5$.}
\end{center}
\end{figure}

\subsection{Cross-sections} We say that a set of sentences
$A$ is a {\em cross-section} of a set of sentences $S$ if 
$A\subset S$ and for each
sentence
$s\in S$ there exists exactly one sentence in $A$
equivalent to $s$.  All the cross-sections of $S$ arise
by a process of selecting exactly one element from each $\sim$-equivalence
class in $S$.

\subsection{Enumeration of Wigner words by Wigner words}  
Fix a
letter
$\alpha$. For each positive integer
$i$ choose a cross-section
$W_i$ of the set of Wigner words so as to
achieve the following conditions:
\begin{itemize}
\item For all $i$, the letter $\alpha$ appears in no word
belonging to $W_i$.
\item For all distinct $i$ and $j$, every word belonging to $W_i$ is
disjoint from every word belonging to $W_j$.
\end{itemize}
(This is always possible to achieve
because the set of letters is countably infinite.)
Let $\varphi$ be any real-valued function of words 
such that $\varphi(w)$ depends
only on the equivalence class of $w$ and vanishes for $\ell(w)\gg 0$,
in which case the support of $\varphi$ consists of only finitely
many equivalence classes of words.
By Proposition~\ref{Proposition:ExtremalWeakWigner}(iii) we have an
enumeration formula
\begin{equation}\label{equation:WignerAutoEnumeration}
\sum_{w}\varphi(w)=
\sum_{r=0}^\infty 
\sum_{w_1\in W_1}\cdots \sum_{w_r\in W_r}\varphi(\alpha w_1\cdots \alpha
w_r\alpha)
\end{equation}
where $w$ ranges over any cross-section of the set of Wigner
words. Formula
(\ref{equation:WignerAutoEnumeration}) leads to many useful recursions. For
example, it implies that the number of equivalence classes of Wigner words
of length $2n+1$ is the $n^{th}$ {\em Catalan number}
$\frac{1}{n+1}\left(\begin{subarray}{c}
2n\\n\end{subarray}\right)$. The latter fact is anyhow clear from the
rooted planar tree interpretation of equivalence classes of Wigner words.

\begin{Proposition}\label{Proposition:CriticalMomentOkay}
Let $w$ be a critical weak Wigner word. Put
$$G=(V,E)=G_w=(V_w,E_w).$$
The following hold:
\begin{enumerate}
\item $G$ is connected.
\item Either $\card V-1=\card E$ or $\card V=\card E$.
\item If $\card V-1=\card E$, then:
\begin{enumerate}
\item $G$ is a tree.
\item With exactly one exception $w$ visits each
edge of $G$ exactly twice.
\item  But $w$ visits the exceptional edge
exactly four times.
\end{enumerate}
\item If $\card V=\card E$, then:
\begin{enumerate}
\item $G$ is not a tree.
\item  $w$ visits each edge of $G$ exactly
twice. 
\end{enumerate}
\end{enumerate}
\end{Proposition}
\noindent We state these facts for the sake of convenient reference. We
omit the easy proofs.

\begin{Proposition}\label{Proposition:Marriage}
Let $a=[w_i]_{i=1}^n$ be a weak
CLT sentence consisting of $n$ words. 
(i) We have 
$\weight a\leq \sum_{i=1}^n\frac{\ell(w_i)-1}{2}$.
(ii) Suppose now that equality holds, i.~e., that $a$ is a CLT sentence.
Then
the words $w_i$ of the sentence $a$
are perfectly matched in the sense that for all $i$
there exists unique $j$ distinct from $i$ such
that $w_i$ and $w_j$ have a letter in common. 
In particular, $n$ is even. 
\end{Proposition}
\noindent This assertion (without proof) was made in
\cite{jonsson}.

\proof  Lemma~\ref{Lemma:Marriage} below
is the essential point of the proof.
\qed

\begin{Lemma}\label{Lemma:Marriage}
Let 
$a=[w_i]_{i=1}^n$
be a weak
CLT sentence consisting of $n$ words. Put $G=G_a$.
Let $k$ be the
number of connected components of $G$. Then
(i) $k\leq\left\lfloor \frac{n}{2}\right\rfloor$
and (ii) $\weight a\leq
k-n+\left\lfloor\frac{\sum_{i=1}^n\ell(w_i)}{2}\right\rfloor$,
where $\lfloor x\rfloor$ denotes the greatest integer less than or
equal to $x$.
\end{Lemma}
\proof 
Inequality (i) is trivial: by hypothesis every word $w_i$ of the sentence
$a$ is ``mated'' with at least one other word $w_j$ ($j\neq i$) of the
sentence in the sense that the connected subgraphs $G_{w_i}$ and $G_{w_j}$
share an edge and {\em a fortiori} share a vertex.

Harder work is required to prove inequality (ii).
Put $a=[[\alpha_{ij}]_{j=1}^{\ell(w_i)}]_{i=1}^n$,  $I=\bigcup_{i=1}^n
\{i\}\times \{1,\dots,\ell(w_i)-1\}$ and
$A=[\{\alpha_{ij},\alpha_{i,j+1}\}]_{(i,j)\in I}$. We visualize $A$ as a
left-justified table of $n$  rows. 
 Let $G'=(V',E')$ be any spanning forest in $G$. 
Since every
connected component of $G'$ is a tree, we have
$\weight a=k+\card E'$ and so in order to prove 
(ii), we just have to bound $\card E'$.
Now let $X=\{X_{ij}\}_{(i,j)\in I}$ be a table of the same ``shape'' as
$A$, but with all entries equal either to $0$ or $1$.  We call
$X$ an {\em edge-bounding table} under the following conditions:
\begin{itemize}
\item For all $(i,j)\in I$, if $X_{ij}=1$, then
$A_{ij}\in E'$. 
\item For each $e\in E'$ there exist distinct 
$(i_1,j_1),(i_2,j_2)\in I$ such that  $X_{i_1j_1}=X_{i_2j_2}=1$
and $A_{i_1j_1}=A_{i_2j_2}=e$.
\item For each $e\in E'$ and index $i\in \{1,\dots,n\}$,
if $e$ appears in the $i^{th}$ row of $A$
then there exists $(i,j)\in I$ such that
$A_{ij}=e$ and $X_{ij}=1$.
\end{itemize}
For any edge-bounding table $X$ the corresponding quantity
$\frac{1}{2}\sum_{(i,j)\in I} X_{ij}$ bounds $\card E'$, whence the
terminology. At least one edge-bounding table exists, namely the table
with a
$1$ in position
$(i,j)$  for each $(i,j)\in I$ such that
$A_{ij}\in E'$ and $0$'s elsewhere.
Now let $X$ be an edge-bounding table such
that for some index $i_0$ all the entries of $X$ in the $i_0^{th}$
row are equal to $1$. Then the closed word $w_{i_0}$ is a walk in $G'$,
and hence
by the parity principle  every entry in the $i_0^{th}$ row
of $A$ appears there an even number of times and {\em a fortiori} at least
twice. Now choose
$(i_0,j_0)\in I$ such that  $A_{i_0j_0}\in E'$ appears in more than one
row of $A$. Let $Y$ be the table obtained by replacing the entry $1$ of $X$
in position $(i_0,j_0)$ by the entry $0$.
Then it is not difficult to check that
 $Y$ is again an edge-bounding table. Proceeding in this way we can
find an edge-bounding table with $0$ appearing at least
once in every row, and hence we have  $\card E'\leq \lfloor
\frac{\card I-n}{2}\rfloor$, which is exactly what we need to prove (ii).
\qed

\subsection{Enumeration of CLT sentences by CLT word-pairs}
Fix an even positive integer $n$ and
for
$i=1,\dots,n/2$ choose a cross-section $P_i$ of the set of CLT word-pairs
so
as to achieve the following condition:
\begin{itemize}
\item For all distinct $i$ and $j$, every word-pair belonging to $P_i$ is
disjoint from every word-pair belonging to $P_j$.
\end{itemize}
We declare a permutation $\sigma$ of $\{1,\dots,n\}$ to be a {\em
perfect matching} if it satisfies the following conditions:
\begin{itemize}
\item $\sigma(2i-1)<\sigma(2i)$ for $i=1,\dots,n/2$.
\item $\sigma(2i-1)<\sigma(2i+1)$ for $i=1,\dots,n/2-1$.
\end{itemize}
Now let
$\varphi$ be any real-valued function of $n$-word-long sentences
$a=[w_i]_{i=1}^{n}$ such that $\varphi(a)$
depends only on the equivalence class of
$a$ and vanishes for $\sum_{i=1}^{n} \ell(w_i)\gg 0$,
in which case the support of $\varphi$ consists of only finitely
many equivalence classes of sentences.
By Proposition~\ref{Proposition:Marriage} we have an enumeration formula
\begin{equation}\label{equation:Marriage}
\sum_{a} \varphi(a)=
\sum_{[p_1,p_2]\in P_1}\cdots \sum_{[p_{n-1},p_{n}]\in
P_{n/2}}\sum_{\begin{subarray}{c}
\sigma\in S_{n}\\
\mbox{\tiny $\sigma$: perfect matching}
\end{subarray}}
\varphi([p_{\sigma^{-1}(i)}]_{i=1}^{n})
\end{equation}
where $a$ ranges over any cross-section of the set of $n$-word-long CLT
sentences.

\begin{Proposition}\label{Proposition:CLTWordPairMomentControl}
Let $a=[w,x]$ be  a CLT word-pair and put 
$$G=(V,E)=G_{a}=(V_a,E_a).$$
For each $e\in E$ let $\nu(e,w)$ (resp., $\nu(e,x)$) denote the
number of visits to $e$ by the word/walk $w$ (resp., $x$).
The following hold:
\begin{enumerate}
\item $G$ is connected. 
\item Either $\card V-1=\card E$ or $\card V=\card E$.
\item If $\card V-1=\card E$, then:
\begin{enumerate}
\item $G$ is a tree.
\item For all $e\in E$ both $\nu(e,w)$ and $\nu(e,x)$ are even.
\item For unique $e_0\in E$ we have $\nu(e_0,w)=\nu(e_0,x)=2$.
\item For all $e\in E\setminus\{e_0\}$ we have
$\nu(e,w)+\nu(e,x)=2$.
\item Both $w$ and $x$ are Wigner words.
\end{enumerate}
\item If $\card V=\card E$, then:
\begin{enumerate} 
\item $G$ is not a tree.
\item For all $e\in E$ we have $\nu(e,w)+\nu(e,x)=2$.
\item For some $e\in E$ we have $\nu(e,w)=\nu(e,x)=1$. 
\end{enumerate}
\end{enumerate}
\end{Proposition}
\noindent We state these facts for the sake of convenient reference. We
omit the easy proofs. 

\section{Limit calculations}
\label{section:Limcal}
We work out  limits of and estimates for moments needed as ``raw
material'' for the proofs of Theorems~\ref{Theorem:Preliminary}
and \ref{Theorem:Main}. Assumption~\ref{ass-main} remains in force
throughout these calculations.

\subsection{Random variables indexed by sentences}
\label{subsection:SentenceIndexed} 
Fix a sentence 
$$a=[w_i]_{i=1}^n=[[\alpha_{ij}]_{j=1}^{\ell(w_i)}]_{i=1}^n$$
consisting of $n$ words. We attach several random variables to $a$, 
as follows.
\subsubsection{}
We define
$$\xi(a)=\prod_{i=1}^n \prod_{j=1}^{\ell(w_i)-1}
\xi_{\{\alpha_{ij},\alpha_{i,j+1}\}}.$$ From the independence of the
family
$\{\xi_e\}$ and assumption (\ref{equation:MomentConstraint}) concerning
the absolute moments of these random variables, we deduce that
\begin{equation}\label{equation:FundamentalBound}
E|\xi(a)|=\prod_{
e:\,\mbox{\tiny
edge of
$G_a$}}E|\xi_e|^{\nu(e)}\leq\prod_{\begin{subarray}{c}
e=\{\alpha,\beta\},\\
\mbox{\tiny
edge of
$G_a$}
\end{subarray}}
\left\{\begin{array}{cl}
s^{(\nu(e))}(\kappa_0(\alpha),\kappa_0(\beta))&\mbox{if
$\alpha\neq \beta$,}\\ d^{(\nu(e))}(\kappa_0(\alpha))&\mbox{if
$\alpha=\beta$,}
\end{array}\right.
\end{equation}
where  
\begin{equation}
\label{eq:040504a}
\nu(e)=\left(\begin{array}{l}
\mbox{total
number of visits made to $e$ }\\
\mbox{by all the words/walks $w_i$}
\end{array}\right).
\end{equation}
(While $\nu(e)$ depends on the sentence $a$, to avoid unnecessary clutter
we omit this dependence 
from the notation). 
Further and crucially, since all the random variables of the family
$\{\xi_e\}$ are of mean zero, if
$w$ is a closed word, then
$E\xi(w)=0$ unless
$w$ is a weak Wigner word.

\subsubsection{}
We define
$$\bar{\xi}(a)=\prod_{i=1}^n (\xi(w_i)-E\xi(w_i)).$$
Expanding the product on the right in evident fashion we find that
\begin{equation}\label{equation:XiBarExpansion}
\bar{\xi}(a)=
\sum_{
I\subset\{1,\dots,n\}}(-1)^{\card I}
\prod_{j\in \{1,\dots,n\}\setminus
I}\xi(w_j)\cdot \prod_{i\in I}E\xi(w_i).
\end{equation}
Clearly $E|\bar{\xi}(a)|$ is
bounded by a constant depending only on $\sum_{i=1}^n \ell(w_i)$.
Further and crucially, if all the words $w_i$ are closed, then
we have
$E\bar{\xi}(a)=0$ unless $a$ is a weak CLT sentence.

\subsubsection{Auxiliary color-valued random variables}
We fix a letter-indexed i.i.d.\ family  $\{\kappa(\alpha)\}$ of
color-valued random variables with common distribution $\theta$.
These random variables are going to be used only for bookkeeping purposes.
They need not be defined on the same probability space as
the random variables $\xi_e$.

\subsubsection{}
\label{sec-almost0}
Put
$$
M(a)
=\prod_{\begin{subarray}{c}
e=\{\alpha,\beta\},\\
\mbox{\tiny
edge of
$G_a$}
\end{subarray}}
\left\{\begin{array}{cl}
0&\mbox{if $\nu(e)=1$,}\\
s^{(\nu(e))}(\kappa(\alpha),
\kappa(\beta))&\mbox{if $\nu(e)>1$ and $\alpha\neq \beta$,}\\
d^{(\nu(e))}(\kappa(\alpha))
&\mbox{if $\nu(e)>1$ and $\alpha=\beta$,}
\end{array}\right.
$$ 
where $\nu(e)$ is as in (\ref{eq:040504a}). 
\subsubsection{} Put
$$\overline{M}(a)=
\sum_{
I\subset\{1,\dots,n\}}
(-1)^{\card I}M(a/I)\cdot
\prod_{i\in I} M(w_i),
$$
where for $I\neq \{1,\dots,n\}$ we denote by $a/I$ the sentence obtained
by striking the $i^{th}$ word of $a$ for all $i\in I$, and 
for $I=\{1,\dots,n\}$ we agree to put $M(a/I)=1$.
Note the analogy with expansion
(\ref{equation:XiBarExpansion}). 
Note also that $M(a)$, $\overline{M}(a)$ are random variables.

\subsubsection{}\label{sec-almost1}
For each
$n$-tuple
$p=[p_i]_{i=1}^n$ of nonnegative integers put
$$H_p(a)=\sum_\pi\prod_{i=1}^n\prod_{j=1}^{\ell(w_i)}D(\kappa(\alpha_{ij}))
^{\pi_{ij}}$$
where $\pi=[[\pi_{ij}]_{j=1}^{\ell(w_i)}]_{i=1}^n$ ranges over families
of nonnegative integers subject to the constraints that
$\sum_{j=1}^{\ell(w_i)}\pi_{ij}=p_i$ for $i=1,\dots,n$. 
Note that
$H_p(a)=\prod_{i=1}^nH_{p_i}(w_i)$. It is convenient to set $H_p(a)=0$ for
every $n$-tuple $p$
of integers
such that $p_i<0$ for some $i$. 
We  write
$$MH_p(a)=M(a)H_p(a),\;\;\;\overline{M}H_p(a)=\overline{M}(a)H_p(a)$$
in order to abbreviate notation. 

\subsubsection{}\label{subsubsection:FactorizationAndIndependence}
Note that $MH_p(a)$ (resp., $\overline{M}H_p(a)$) remains
unchanged if for some permutation
$\sigma$ of $\{1,\dots,n\}$ we replace $a$ by the sentence
$[w_{\sigma(i)}]_{i=1}^n$ and $p$ by the $n$-tuple
$[p_{\sigma(i)}]_{i=1}^n$.  Note further that if
the sentence $a$ can be presented as the concatenation of pairwise disjoint
sentences
$b_1,\dots,b_k$ where $b_i$ is $n_i$ words long, and correspondingly we
present the $n$-tuple $p$ as the concatenation of tuples $q_1,\dots,q_k$
where
$q_i$ is an
$n_i$-tuple, then
$MH_p(a)=\prod_{i=1}^k MH_{q_i}(b_i)$ (resp.,
$\overline{M}H_p(a)=\prod_{i=1}^k \overline{M}H_{q_i}(b_i)$), and moreover
the factors on the right are independent.

\subsection{Admissibility} 
Let $a=[w_i]_{i=1}^n$ be a sentence consisting of $n$ words.
For each edge $e$ of the graph $G_a$, let $\nu(e)$ be the total
number of visits to $e$ by the words/walks $w_i$. 
We say that $a$ is {\em weakly admissible} if 
for all edges $e$ of $G_a$ the following hold:
\begin{itemize}
\item $\nu(e)\in \{1,2,4\}$.
\item If $\nu(e)=4$,
then $e$ is nondegenerate.
\end{itemize}
We say that $a$ is {\em admissible} if for every nonempty
subset $\{i_1<\dots<i_\ell\}\subset\{1,\dots,n\}$
the subsentence $[w_{i_\nu}]_{\nu=1}^\ell$ is weakly admissible.
For words weak admissibility and admissibility are the same
thing. By Proposition~\ref{Proposition:ExtremalWeakWigner} every
Wigner word is admissible. By
Proposition~\ref{Proposition:CriticalMomentOkay} every critical weak
Wigner word is admissible. By Propositions~\ref{Proposition:Marriage}
and \ref{Proposition:CLTWordPairMomentControl} every CLT sentence
is admissible.

\begin{Proposition}\label{Proposition:Randomization}
Let $a=[w_i]_{i=1}^n$ be a weakly admissible sentence consisting of $n$
words.  Let $p=[p_i]_{i=1}^n$ be an $n$-tuple of nonnegative integers.
Let $\gamma_1,\dots,\gamma_r$ be distinct letters such that $\supp a
\subset\{\gamma_1,\dots,\gamma_r\}$.
Then there exists
a function
$f$ on the product of
$r$ copies of color space with the following properties:
\begin{itemize}
\item $f$ is bounded and measurable.
\item $f$ has discontinuity set of measure
zero with respect to $\theta^{\otimes r}$.
\item $f(\kappa(\gamma_1),\dots,\kappa(\gamma_r))=MH_p(a)$.
\item For all distinct letters $\delta_1,\dots,\delta_r$,
the equivalent word $b=[[\beta_{ij}]_{j=1}^{\ell(w_i)}]_{i=1}^{n}$
to which
$a$ codes by the rule
$\gamma_i\mapsto \delta_i$ for $i=1,\dots,r$ satisfies the equation
$$f(\kappa_0(\delta_1),\dots,\kappa_0(\delta_r))=E\xi(b)
\sum_\pi
\prod_{i=1}^n\prod_{j=1}^{\ell(w_i)}D(\kappa_0(\beta_{ij}))^{\pi_{ij}}$$
where
$\pi=[[\pi_{ij}]_{j=1}^{\ell(w_i)}]_{i=1}^{n}$ ranges over
families of nonnegative integers  subject to the constraints
$\sum_{j=1}^{\ell(w_i)}\pi_{ij}=p_i$ for $i=1,\dots,n$.
\end{itemize}
\end{Proposition}
\proof For simplicity we discuss only the case $p=0$ and leave the
remaining details to the reader. Put $G=(V,E)=G_a=(V_a,E_a)$
and as above, for all $e\in E$, let $\nu(e)$ be the total number of visits
to
$e$ made by the words/walks $w_i$. Put
$$
f(c_1,\dots,c_r)=\prod_{e=\{\alpha,\beta\}\in E}\;
\left\{\begin{array}{cl}0&\mbox{if $\nu(e)=1$,}\\
s^{(\nu(e))}(c_{\gamma^{-1}(\alpha)},c_{\gamma^{-1}(\beta)})&\mbox{if
$\nu(e)>1$ and $\alpha\neq \beta$,}\\
d^{(\nu(e))}(c_{\gamma^{-1}(\alpha)})&\mbox{if $\nu(e)>1$
and $\alpha=\beta$,}
\end{array}\right.
$$
where $\gamma^{-1}$ is the inverse of the bijection $(i\mapsto
\gamma_i):\{1,\dots,r\}\rightarrow\{\gamma_1,\dots,\gamma_r\}$. Clearly
$f$ has the first three of the desired properties. If $\nu(e)=1$ for some
$e\in E$, then the fourth property holds trivially (both sides of
the desired
equation vanish identically). Otherwise, if $\nu(e)>1$ for all $e\in E$,
then $f$ has the
fourth property because, under the hypothesis of weak admissibility, we are
operating in the regime in which we enforce equality in the moment bound
(\ref{equation:MomentConstraint}).
\qed

\subsection{Limiting behavior of 
$\langle \overline{L}(\NNN),x^n\rangle$}
Fix a  nonempty finite
 set $\NNN$ of letters and a positive integer $n$.
\subsubsection{}
 We
have an expansion \\
$$\trace X(\NNN)^n=
\sum_w\sum_J
\sum_x N^{\frac{1}{2}(\card J-n)}
\xi(x/J)\prod_{j\in J}D(\kappa_0(\beta_j))
$$
where:\\
\begin{itemize}
\item $w=[\alpha_i]_{i=1}^{n+1}$ ranges over a cross-section of the set of
closed words of length $n+1$;
\item $J$ ranges over subsets of the set
$\{j\in\{1,\dots,n\}\vert\alpha_j=\alpha_{j+1}\}$; 
\item $x=[\beta_i]_{i=1}^{n+1}$ ranges over words 
such that $x\sim w$ and $\supp x\subset \NNN$; and
\item $x/J$ denotes the word obtained by striking the $j^{th}$
letter of $x$ for each $j\in J$.
\\
\end{itemize}
Note that $x/J$ arises from $x$ by selective suppression
of repeated letters. Note also that $E\xi(x/J)=0$ unless $x/J$ is a weak
Wigner word. By considering how we may insert repetitions of letters into a
given weak Wigner word, and after some further algebraic manipulation,
we obtain an expansion
\begin{eqnarray}
\label{equation:WordExpansion2}
\langle \overline{L}(\NNN),x^n\rangle&=&
\sum_w N^{\weight w-\frac{1+\ell(w)}{2}}
\sum_x
\sum_\pi  N^{-\weight
w}E\xi(x)\prod_{i=1}^{\ell(w)}D(\kappa_0(\beta_i))^{\pi_i}
\nonumber\\
&=:&
\sum_w N^{\weight w-\frac{1+\ell(w)}{2}}
S(\NNN,w)
\end{eqnarray}
where:
\begin{itemize}
\item $w$ ranges over a cross-section of the set of weak Wigner words of
length $\leq n+1$;
\item $x=[\beta_i]_{i=1}^{\ell(w)}$ ranges over words such that
$x\sim w$ and $\supp x\subset \NNN$; 
\item $\pi=[\pi_i]_{i=1}^{\ell(w)}$ ranges over $\ell(w)$-tuples
of nonnegative integers summing to\\  $n+1-\ell(w)$; and
\item $S(\NNN,w)$ is the result of carrying out the inner summations
on $x$ and $\pi$.
\end{itemize}
Note that for $n$ fixed, 
as $N\rightarrow \infty$
and $\theta_{\NNN}\Rightarrow \theta$, only the part of the sum
indexed by Wigner words $w$ contributes nonnegligibly.
\subsubsection{}\label{subsubsection:ApproximateIntegrationExample}
In this paragraph fix
attention on a Wigner word $w$ such that
$\ell(w)\leq n+1$.
We want
to understand the subsum $S(\NNN,w)$ appearing in formula
(\ref{equation:WordExpansion2}) as a function of $\NNN$. Let
$\gamma_1,\dots,\gamma_r$ be an enumeration  of
$\supp w$. Since Wigner words are admissible,
Proposition~\ref{Proposition:Randomization} provides us with a function
$f$ defined on the product of $r$ copies of color space with the following
properties:\\
\begin{itemize}
\item $f$ is bounded and measurable.\\
\item $f$ has discontinuity set of measure
zero with respect to $\theta^{\otimes r}$.\\
\item $MH_{n+1-\ell(w)}(w)=f(\kappa(\gamma_1),\dots,\kappa(\gamma_r))$.\\
\item $\displaystyle S(\NNN,w)=N^{-r}\sum_{\begin{subarray}{c}
(\beta_1,\dots,\beta_r)\in\NNN^r\\
\mbox{\tiny $\beta_1,\dots,\beta_r$: distinct}
\end{subarray}}f(\kappa_0(\beta_1),\dots,\kappa_0(\beta_r))$.\\
\end{itemize} Now let
$[\NNN_k]_{k=1}^\infty$ be as in Assumption
\ref{ass-main}. We clearly have
$$\lim_{k\rightarrow\infty}S(\NNN_k,w)=
\int\cdots \int f(c_1,\dots,c_r)\theta(dc_1)\cdots
\theta(dc_r)=EMH_{n+1-\ell(w)}(w).
$$
We remark that it is here we make use of the hypothesis that color space
is Polish: we need it to guarantee weak convergence
$\theta_{\NNN_k}^{\otimes r}\Rightarrow
\theta^{\otimes r}$.
\subsubsection{}
We may now conclude
that
\begin{equation}\label{equation:WordExpansion3}
\lim_{k\rightarrow\infty}
\langle \overline{L}(\NNN_k),x^n\rangle=
\sum_w EMH_{n+1-\ell(w)}(w)
\end{equation}
where the sum on
the right is extended over a cross-section of the set of Wigner words.
Note that only finitely many terms on the right are nonvanishing  because
$p<0\Rightarrow H_p\equiv 0$.

\begin{Lemma}\label{Lemma:WeakOverlineConvergence}
With $C$ as in (\ref{equation:MuSupportBound}),
$\overline{L}(\NNN_k)$ converges weakly to a limit $\mu$ 
supported in the interval $[-C,C]$,
and moreover
$\langle\overline{L}(\NNN_k),x^n\rangle\rightarrow
\langle \mu,x^n\rangle$ for all integers $n>0$.
\end{Lemma}
\proof It is enough to prove  that the right side of
(\ref{equation:WordExpansion3}) is $O(C^n)$.
There are 
$\left(\begin{subarray}{c}p+n-1\\
n-1\end{subarray}\right)$
$n$-tuples of nonnegative integers
summing to $p$. Consequently
we have
$$|MH_{p}(w)|\leq
\left(\begin{subarray}{c}p+\ell(w)-1\\\ell(w)-1\end{subarray}\right)
(|s^{(2)}|_\infty^{1/2})^{\ell(w)-1}|D|_\infty^{p}
$$
for all Wigner words $w$ and nonnegative integers $p$. 
There are
$\frac{1}{\ell+1}\left(\begin{subarray}{c}
2\ell\\\ell\end{subarray}\right)$ equivalence classes of Wigner words
of length $2\ell+1$ and clearly there are no Wigner words of
even length. Consequently  there are
$O(2^n)$ equivalence classes of Wigner words of length $\leq
n+1$. The desired $O(C^n)$ bound for the right side of
(\ref{equation:WordExpansion3}) follows.
\qed

\subsection{Limiting behavior of
$E\prod_{i=1}^n (\trace X(\NNN)^{\nu_i}-E\trace
X(\NNN)^{\nu_i})$}
Again fix a finite non-empty set $\NNN$ of
letters and a positive integer $n$. Also fix positive integers
$\nu_1,\dots,\nu_n$ and put $\nu=[\nu_i]_{i=1}^n$.
\subsubsection{}
 We have an expansion
$$\begin{array}{cl}
&\displaystyle\prod_{i=1}^n (\trace X(\NNN)^{\nu_i}-E\trace
X(\NNN)^{\nu_i})\\\\
=&\displaystyle
\sum_a\sum_K\sum_b
 N^{\frac{1}{2}\sum_{i=1}^n
(\card K_i-\nu_i)}\prod_{i=1}^n\left(
\bar{\xi}(x_i/K_i)\prod_{j\in K_i}D(\kappa_0(\beta_{ij}))\right)
\end{array}
$$
where:\\
\begin{itemize}
\item $a=[w_i]_{i=1}^n=[[\alpha_{ij}]_{j=1}^{\ell(w_i)}]_{i=1}^n$ ranges
over a cross-section of the set of sentences $n$ words long 
with $i^{th}$ word of length $\nu_i+1$ for $i=1,\dots,n$;
\item $K=[K_i]_{i=1}^n$
 ranges over $n$-tuples of sets of positive
integers such that $K_i$ is a subset of
$\{j\in \{1,\dots,\nu_i\}\vert
\alpha_{ij}=\alpha_{i,j+1}\}$ for
$i=1,\dots,n$; 
\item $b=[x_i]_{i=1}^n=[[\beta_{ij}]_{j=1}^{\nu_i+1}]_{i=1}^n$ ranges
over  sentences $b\sim a$ such that $\supp b\subset \NNN$; and
\item $x_i/K_i$ denotes the word obtained by striking the $k^{th}$
letter of $x_i$ for all $k\in K_i$. \\
\end{itemize}
After some further algebraic
manipulation, we obtain an
expansion
\begin{equation}\label{equation:SentenceExpansion2}
\begin{array}{cl}
&\displaystyle E\prod_{i=1}^n (\trace X(\NNN)^{\nu_i}-E\trace
X(\NNN)^{\nu_i})\\\\ =&\displaystyle\sum_a
 N^{\weight a-\sum_{i=1}^n\frac{\ell(w_i)-1}{2}
}\sum_b 
\sum_\pi N^{-\weight a}E\bar{\xi}(b)
\prod_{i=1}^n\prod_{j=1}^{\ell(w_i)}D(\kappa_0(\beta_{ij}))^{\pi_{ij}}
\end{array}
\end{equation}
where:
\begin{itemize}
\item $a=[w_i]_{i=1}^n$ ranges over a cross-section of the set of
weak CLT sentences $n$ words long with $i^{th}$ word of length $\leq
\nu_i+1$ for
$i=1,\dots,n$;
\item $b=[[\beta_{ij}]_{j=1}^{\ell(w_i)}]_{i=1}^n$ ranges over sentences
$b\sim a$ such that $\supp b\subset \NNN$;  and
\item $\pi=[[\pi_{ij}]_{j=1}^{\ell(w_i)}]_{i=1}^n$ ranges
over families of nonnegative integers subject to the constraints
$\sum_{j=1}^{\ell(w_i)}\pi_{ij}=\nu_i+1-\ell(w_i)$ for 
$i=1,\dots,n$.
\end{itemize}
Note that for fixed $\nu$,
 as $N\rightarrow \infty$
and $\theta_{\NNN}\Rightarrow \theta$, only the part of the sum
indexed by CLT sentences $a$ contributes nonnegligibly.
\subsubsection{}

Now let $[\NNN_k]_{k=1}^\infty$ be as in Assumption \ref{ass-main}.
Since CLT words are admissible,
an analysis similar
to that undertaken in
\S\ref{subsubsection:ApproximateIntegrationExample} leads to the
conclusion that
\begin{equation}
\label{equation:SentenceExpansion3}
\lim_{k\rightarrow\infty}
E\prod_{i=1}^n (\trace X(\NNN_k)^{\nu_i}-E\trace
X(\NNN_k)^{\nu_i})= \sum_a
E\overline{M}H_{[\nu_i+1-\ell(w_i)]_{i=1}^n}(a)
\end{equation}
where $a=[w_i]_{i=1}^n$ ranges over a cross-section of the set of CLT
sentences
$n$ words long. Since the analysis is straightforward, somewhat long, and
very tedious, we omit it. Note that only finitely
many nonzero terms appear in the sum on the right.

\begin{Lemma}\label{Lemma:GaussMarriage}
There exists a family $[Y_n]_{n=1}^\infty$ of mean zero random
variables defined on a common probability space with Gaussian joint
distribution such that for all positive integers $n$ and positive
integers
$\nu_1,\dots,\nu_n$  the right side of limit
formula (\ref{equation:SentenceExpansion3}) gives the expectation
$E\prod_{i=1}^nY_{\nu_i}$.
\end{Lemma}
\proof Let $A(\nu_1,\dots,\nu_n)$ denote the right side 
of (\ref{equation:SentenceExpansion3}).
The matrix $[[A(i,j)]_{i=1}^\infty]_{j=1}^\infty$
is symmetric and every finite block $[[A(i,j)]_{i=1}^r]_{j=1}^r$
in the upper left corner
is positive semidefinite since it is the limit of such matrices.
Consequently there exists a family $[Y_n]_{n=1}^\infty$ of mean zero 
random variables on a common probability space
with Gaussian joint distribution such that 
$EY_iY_j=A(i,j)$ for all $i$ and $j$.
By the enumeration formula (\ref{equation:Marriage}) 
and the relations discussed in
\S\ref{subsubsection:FactorizationAndIndependence}, we have
$$
A(\nu_1,\dots,\nu_n)= 
\left\{\begin{array}{cl}
\displaystyle
\sum_{\begin{subarray}{c}
\sigma\in
S_n\\
\mbox{\tiny
$\sigma$:
perfect
matching}
\end{subarray}}\prod_{i=1}^{n/2}A(\nu_{\sigma(2i-1)},\nu_{\sigma(2i)})
&\mbox{if
$n$ is even,}\\\\ 0&\mbox{if $n$ is odd.}\\
\end{array}\right.
$$
But the expression on the right side is the Wick
formula for the expectation
$E\prod_{i=1}^n Y_{\nu_i}$, cf.~\cite[Theorem 1.28]{janson}. 
\qed

\section{Proofs of Lemma \ref{LemmaTheorem:Preliminary}
and Theorem~\ref{Theorem:Preliminary}}
\label{sec:proofI}

\begin{Lemma}\label{Lemma:CauchyTwiceChebyshevOnce}
Fix
$K>\max(1,C^2)$ with $C$ as in (\ref{equation:MuSupportBound}). Then we
have
$$\lim_{k\rightarrow\infty} \langle
\overline{L}(\NNN_k),|g|\one_{|x|>K}\rangle=0$$
for every real-valued measurable function $g$ on the
real line with polynomial growth at infinity.
\end{Lemma}
\proof Let $n$ be any nonnegative integer. We have
$$
\langle \overline{L}(\NNN_k),|x|^n\one_{|x|>K}\rangle
\leq \sqrt{\langle \overline{L}(\NNN_k),x^{2n}\rangle}
\sqrt{\langle \overline{L}(\NNN_k),\one_{|x|>K}\rangle} 
\leq \langle \overline{L}(\NNN_k),x^{2n}\rangle/K^n
$$
by Cauchy-Schwartz followed by Chebyshev, and hence
\begin{equation}\label{equation:OfersEstimate}
\limsup_{k\rightarrow\infty}\langle
\overline{L}(\NNN_k),|x|^n\one_{|x|>K}\rangle
\leq \limsup_{k\rightarrow\infty}
\langle \overline{L}(\NNN_k),x^{2n}\rangle/K^n.
\end{equation}
Because $K>1$, the quantity on the left side of
(\ref{equation:OfersEstimate}) is an increasing function of $n$, and
moreover 
that quantity bounds
$\limsup_{k\rightarrow\infty}\langle
\overline{L}(\NNN_k),|g|\one_{|x|>K}\rangle$
for all $n\gg 0$ (because $g$ is of polynomial growth). But by 
Lemma~\ref{Lemma:WeakOverlineConvergence}, because $K>C^2$, the
quantity on the right side of (\ref{equation:OfersEstimate})  tends to
$0$ as
$n\rightarrow\infty$. The result follows. \qed

\subsection{The functions $\Phi^{(w,p)}(c)$}
To each Wigner word $w$ and nonnegative integer $p$ we
associate a real-valued bounded measurable function $\Phi^{(w,p)}(c)$ on
color space by the following recursive procedure. As in
Proposition~\ref{Proposition:ExtremalWeakWigner}, in the unique way
possible, write 
$w=\alpha w_1\cdots
\alpha w_r\alpha$ where
$\alpha$ is the first letter of $w$ and the $w_i$ are pairwise disjoint
Wigner words in which
$\alpha$ does not appear, and then put
\begin{equation}\label{equation:PhiAuxiliary}
\Phi^{(w,p)}(c)=\sum_\pi D(c)^{\pi_{0}+\cdots +\pi_{r}}\prod_{i=1}^r
\int s^{(2)}(c,c')\Phi^{(w_{i},\pi_{i+r})}(c')\theta(dc')
\end{equation}
where $\pi=[\pi_i]_{i=0}^{2r}$ ranges over $(2r+1)$-tuples of
nonnegative integers summing to $p$. 
By convention, if $w$ is the single letter word $\alpha$, then
$r=0$ and therefore $\Phi^{(\alpha,p)}=D(c)^p$, which gives
a way to initialize the recursions (\ref{equation:PhiAuxiliary}).
Note that for fixed $p$ and
$c$ the quantity
$\Phi^{(w,p)}(c)$ depends only on the equivalence class of $w$.
Intuitively, $\Phi^{(w,p)}(c)$ determines the dominant 
contribution
to the expectation of $\trace X(\NNN)^{\ell(w)+p}$
by those terms that use entries from $D$ $p$ times, such that
when these are discarded, the resulting word determined by the
indices is
equivalent to $w$, and such that the color
of the initial letter is  $c$. 
For example, in the special case that $D(\cdot)=0$, one must have
$p=0$, hence all $\pi_i$ vanish, and  
the contribution, for a given $w$,
can be visualized by writing on each
edge $(v_1,v_2)$ of the rooted  planar tree the value
$(s^{(2)})^{1/2}(\kappa(v_1),
\kappa(v_2))$,
collecting the product of such values along the  exploration
path determined by the word $w$, and averaging over the choices of colors 
{\it except} for the choice of the color of the root, which is fixed at
$c$.

\begin{Lemma}\label{Lemma:PhiDecomposition}
We have the following identity
of formal power series in $t$ with coefficients in the space of
real-valued bounded measurable functions on color space:
\begin{equation}\label{equation:PhiDecomposition}
\Phi(c,t)=\sum_w \sum_{p=0}^\infty \Phi^{(w,p)}(c)t^{\ell(w)+p}
\end{equation}
Here $w$ ranges over a cross-section of the set of Wigner words.
\end{Lemma}
\proof Via the enumeration formula (\ref{equation:WignerAutoEnumeration})
it
follows  from definition
(\ref{equation:PhiAuxiliary}) that the power series on the right side
of (\ref{equation:PhiDecomposition})  satisfies 
(\ref{equation:PhiRecursion}), whence the result.
\qed

\begin{Lemma}\label{Lemma:ConditionalPhiInterpretation}
Let $w$ be a Wigner word. 
Let $\alpha$ be the first letter of $w$. 
Let $p$ be a nonnegative integer. Then we have
\begin{equation}\label{equation:ConditionalPhiInterpretation}
E( MH_p(w)
\vert\kappa(\alpha))=\Phi^{(w,p)}(\kappa(\alpha)), \quad \mbox{\rm a.s.}.
\end{equation}
\end{Lemma}
\proof As in definition
(\ref{equation:PhiAuxiliary}), write
$w=\alpha w_1\cdots \alpha w_r\alpha$ where the
$w_i$ are pairwise disjoint Wigner words in which $\alpha$ does not
occur and let $\alpha_i$ denote the first letter of $w_i$. By
definition of $M(\cdot)$ and $H_p(\cdot)$ we have
\begin{equation}\label{equation:ConditionalPhiInterpretationBis}
MH_p(w)=\sum_\pi D(\kappa(\alpha))^{\pi_{0}+\cdots
+\pi_{r}}\prod_{i=1}^r
s^{(2)}(\kappa(\alpha),\kappa(\alpha_i))MH_{\pi_{i+r}}(w_i)
\end{equation}
where
$\pi=[\pi_i]_{i=0}^{2r}$ ranges over $(2r+1)$-tuples of nonnegative
integers summing to
$p$.  Now take conditional expectations 
on both sides of (\ref{equation:ConditionalPhiInterpretationBis}).
By induction on $\ell(w)$, and the relations of independence built into
the definitions of $M(\cdot)$ and $H_\cdot(\cdot)$, we get
(\ref{equation:ConditionalPhiInterpretation})
after a
routine calculation.
\qed

\subsection{Ends of the proofs}
\label{subsection:TheoremOneProof}

\subsubsection{Proof of Lemma \ref{LemmaTheorem:Preliminary}} 
Uniqueness of a probability measure with moments
(\ref{equation:PhiMomentFormula}) and 
support (\ref{equation:MuSupportBound}) (which is compact) is clear.
Only existence requires proof. After enlarging the originally given model
in evident fashion we may
assume without loss of generality that for every letter there exist
infinitely many letters of the same color. And then we may assume without
loss of generality that
$\sigma=\theta$ because in Assumption \ref{ass-main} we may substitute
$\theta_{\NNN}$ for
$\theta$ without falsifying it. 
Now fix any sequence
$[\NNN_k]_{k=1}^\infty$ as in Assumption \ref{ass-main}.
Let $\mu$ be the weak limit of $\overline{L}(\NNN_k)$ provided by
Lemma~\ref{Lemma:WeakOverlineConvergence}. By the cited lemma,
$\mu$ satisfies the support bound
(\ref{equation:MuSupportBound}). Moreover, by the cited lemma
combined with limit formula (\ref{equation:WordExpansion3}), the
measure
$\mu$ has
moments
\begin{equation}\label{equation:WordExpansion4}
\langle \mu, x^n\rangle=\sum_{w\in W} EMH_{n+1-\ell(w)}(w)
\end{equation}
where $w$ ranges over a cross-section of the set of Wigner words.
By Lemmas~\ref{Lemma:PhiDecomposition} and
\ref{Lemma:ConditionalPhiInterpretation} we can evaluate the
right side of (\ref{equation:WordExpansion4}). We find finally that
moment formula (\ref{equation:PhiMomentFormula}) does indeed hold for
$\mu$.
\qed

\subsubsection{Proof of Theorem \ref{Theorem:Preliminary}}
Fix any real-valued bounded continuous function
$f$ on the real line and $\epsilon>0$. For the convergence
$L(\NNN_k)\Rightarrow \mu$
it is enough to show that
\begin{equation}\label{equation:TheoremOneNuff}
\lim_{k\rightarrow\infty}
P(|\langle L(\NNN_k),f\rangle-\langle
\mu,f\rangle|>\epsilon)=0.
\end{equation} 
Fix $K$ as in Lemma~\ref{Lemma:CauchyTwiceChebyshevOnce}.
By the Weierstrass approximation theorem write
$$f=g+Q,\;\;\;\sup_{|x|\leq
K}|g(x)|<\epsilon/4$$
where $Q$ is a polynomial function.
We have
\begin{eqnarray*}
\langle L(\NNN_k),f\rangle-\langle \mu,f\rangle
&=&\Big[\langle
L(\NNN_k),\one_{|x|\leq K}g\rangle-
\langle\mu,\one_{|x|\leq K}g\rangle\Big]
+\langle L(\NNN_k),\one_{|x|>K}g\rangle\\
&&+\Big[\langle\overline{L}(\NNN_k),Q\rangle-\langle
\mu,Q\rangle\Big]
+\Big[\langle L(\NNN_k),Q)\rangle-\langle \overline{L}(\NNN_k),Q\rangle
\Big]
\end{eqnarray*}
and therefore have 
\begin{eqnarray*}
P(|\langle L(\NNN_k),f\rangle-\langle \mu,f\rangle|>\epsilon)
&\leq &P(\langle L(\NNN_k),\one_{|x|>K}|g|)\rangle >\epsilon/6)\\
&&+
P(|\langle\overline{L}(\NNN_k),Q\rangle-\langle
\mu,Q\rangle|>\epsilon/6)\\ 
&&+P(|\langle L(\NNN_k),Q)\rangle-\langle
\overline{L}(\NNN_k),Q\rangle|>\epsilon/6)\\ &:=&\;\;P_1+P_2+P_3.
\end{eqnarray*}
We have
$P_1\rightarrow 0$ by 
Lemma~\ref{Lemma:CauchyTwiceChebyshevOnce}. We have
$P_2\rightarrow 0$ by Lemma~\ref{Lemma:WeakOverlineConvergence}. 
We have $P_3\rightarrow 0$ by limit formula
(\ref{equation:SentenceExpansion3}).
Therefore (\ref{equation:TheoremOneNuff}) does indeed hold.

We finally turn to proving the convergence
$\mu_{\NNN_k}\Rightarrow \mu$. The proof of Lemma
\ref{LemmaTheorem:Preliminary}
shows that the analogue 
\begin{equation}\label{equation:WordExpansion5}
\langle \mu_{\NNN}, x^n\rangle=\sum_w EM_{\NNN}H_{n+1-\ell(w),\NNN}(w)
\end{equation}
of formula (\ref{equation:WordExpansion4}) holds for any nonempty
finite set of letters
$\NNN$, where the 
random variables
$M_\NNN(w)$ and $H_{p,\NNN}(w)$ are defined
by mimicking the definitions of $M(w)$ and $H_p(w)$, only this
time using a  letter-indexed family $\{\kappa_\NNN(\alpha)\}$
of color-valued family i.i.d.\ random variables with common
law $\theta_\NNN$. Note that $M_\NNN(w)$ and $H_{p,\NNN}(w)$ are
uniformly bounded in $\NNN$. Clearly
for each Wigner word
$w$ and nonnegative integer
$p$ we have convergence in distribution
$M_{\NNN_k}H_{p,\NNN_k}(w)
\rightarrow MH_p(w)$, which extends to the convergence of expectations
by bounded convergence. The sum in (\ref{equation:WordExpansion5}) being
over a finite number of terms, it follows that
$\langle
\mu_{\NNN_k},x^n\rangle\rightarrow
\langle
\mu, x^n\rangle$ for all $n$, and in turn
that $\mu_{\NNN_k}\Rightarrow\mu$ since the measures in play here have
uniformly bounded supports. The proof of
Theorem~\ref{Theorem:Preliminary} is complete. \qed

\section{The F\"{u}redi-Koml\'{o}s circle of ideas}
\label{sec:FK}
In this section, we describe a (rough) technique which allows
us to  bound traces of polynomials of our random 
matrices when the degree of the polynomial is allowed to grow with the
dimension of the matrix. The approach we take is inspired by the work of
 F\"{u}redi and Koml\'{o}s \cite{furedi}. We mention in passing 
that for Wigner matrices all of whose entries have even distributions, much
more detailed information is available in
\cite{soshnikov}.

\subsection{FK sentences} 
Let $a=[w_i]_{i=1}^n$ be a sentence of $n$ words. 
We say that $a$ is an {\em FK sentence} under the following conditions:
\begin{itemize}
\item $G_a$ is a tree.
\item Jointly the words/walks $w_i$ visit no edge of $G_a$ more than
twice.
\item For $i=1,\dots,n-1$, the first letter of $w_{i+1}$
belongs to $\bigcup_{j=1}^{i} \supp w_j$.
\end{itemize}
We say that $a$ is an {\em FK word} if $n=1$. 
Any
word admitting interpretation as a walk on a forest visiting no edge of the
forest more than twice is automatically an FK word.
The constituent words of an FK sentence are FK
words.  If an FK sentence is at least two words long,
then the result
of dropping the last word is again
an FK sentence. If the last word of an FK sentence is at least two
letters long, then the result of dropping the last letter of the last word
is again an FK sentence. 

\subsection{The graph $G^1_a$ associated to a sentence}

Given an $n$-word-long sentence $a=[w_i]_{i=1}^n$,
we define
$G^1_a=(V^1_a,E^1_a)$ to be the subgraph of $G_a=(V_a,E_a)$ with
$V^1_a=V_a$ and $E^1_a$ equal to the set of edges $e\in E_a$
such that the words/walks $w_i$ jointly visit $e$ exactly once.
\begin{figure}[h]
\begin{center}
\includegraphics[scale=.49]{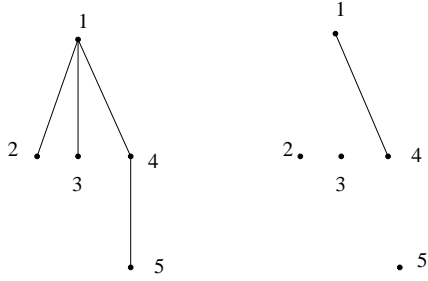}
\caption{The graphs $G_w$ (left) and $G_w^1$ (right) for the FK word
$w=12131454$}
\end{center}
\end{figure}

\begin{Proposition}\label{Proposition:Acronym}
Let $w$ be an FK word.
There is exactly one way to write
$w=w_1\cdots w_r$ where the words $w_i$ are pairwise disjoint Wigner words.
\end{Proposition}
\noindent In this situation, denoting by $\alpha_i$ the first letter of
$w_i$, we declare the word $\alpha_1\cdots\alpha_r$ to be the {\em
acronym} of the FK word $w$. 
\proof
The only possible decomposition $w=w_1\cdots w_r$ of the desired type
is the one with breaks at the edges of $G^1_w$. Since the transition from
$w_{i-1}$ to $w_i$ is along an edge of the tree $G_w$ never again visited
by
$w$, the words
$w_i$ must be pairwise disjoint. Since every edge of $G_w$ visited by
$w_i$ is visited exactly twice by $w$, and the $w_i$ are pairwise
disjoint, in fact 
$w_i$ visits every edge of $G_w$ either twice or never,
hence by the parity
principle
$w_i$ is closed, and hence $w_i$ is a Wigner word.
\qed
\begin{Lemma}\label{Lemma:FKcounting}
There are at most $2^{n-1}$ equivalence classes of FK words
of length $n$.
\end{Lemma}
\proof  From the recursion (\ref{equation:WignerAutoEnumeration})
(see also
\S\ref{subsubsection:PhiChebyshev} below)
it is easy to deduce  that the sum of terms
$t^{\ell(w)}$ extended over a cross-section of the set of Wigner words is 
$$\Phi(t)=\frac{1-\sqrt{1-4t^2}}{2t}.$$ 
Via the preceding lemma, it follows that the sum of terms $t^{\ell(w)}$
extended over a cross-section of the set of FK words is
$$\frac{\Phi(t)}{1-\Phi(t)}=-1/2+\frac{1}{2}\frac{1+2t}{\sqrt{1-4t^2}}
=t+\left(\frac{1}{2}+t\right)\sum_{n=1}^\infty
\left(\begin{array}{c}
2n\\
n
\end{array}\right)t^{2n},
$$
whence the claimed bound.
\qed

\subsection{FK syllabification}
Let $w=[\alpha_i]_{i=1}^n$ be a word of length $n$. Roughly speaking, we
wish to define a
parsing of $w$ into an FK sentence by going sequentially over the letters
in $w$ and  declaring a new word 
each time not doing so would prevent the sentence formed up to that point
from being an FK sentence. More precisely, we
define a sentence $w'$, which we call
the {\em FK syllabification} of $w$, by the following 
procedure. We declare an edge $e$ of $G_w$ to be {\em new} (relative to
$w$) if for some index $1\leq i<n$ we have $e=\{\alpha_i,\alpha_{i+1}\}$
and  $\alpha_{i+1}\not\in \{\alpha_1,\dots,\alpha_i\}$,
and otherwise we declare $e$ to be {\em old}.
We define $w'$ to be the sentence obtained by breaking
$w$ at all visits to old edges of $G_w$
and at third and subsequent visits to new edges of $G_w$.
For example, temporarily spelling with the alphabet
$\{1,2,3\}$, the FK syllabification of
$w=1231$ is the sentence 
$w'=[123,1]$ consisting of two words; the FK syllabification process
has to ``insert a comma'' between $3$ and $1$ because $1231$ is not an FK
word, whereas $1$, $12$ and $123$ are.
It is clear that $G_{w'}$ is a spanning tree in
$G_w$, that  $w'$ is an FK sentence, and that $w$ is the concatenation of
the constituent words of $w'$. Moreover, we have $w=w'$ if and only if $w$
is an FK
word. 
Clearly the FK syllabification process preserves
equivalence, i.~e., $w\sim x\Rightarrow w'\sim x'$.

\begin{Lemma}\label{Lemma:Coding}
Let $a=[w_i]_{i=1}^n$ be a sentence of $n\geq 2$ words.
Put $b=[w_i]_{i=1}^{n-1}$ and $c=w_n$. Assume that
$b$ is an FK sentence, that $c$ is an FK word,
and that the first letter of $c$ belongs to $\supp b$. 
Let
$\gamma_1\cdots \gamma_r$ be the acronym of $c$ spelled out in full.
(Note that by hypothesis $\gamma_1\in \supp b$.)
Let $\ell$ be the largest index such that $\gamma_\ell\in \supp b$
and write $d=\gamma_1\cdots \gamma_\ell$.
The
following conditions are both necessary and sufficient for $a$ to be an
FK sentence: 
\begin{itemize}
\item $d$ is a geodesic in the forest $G^1_b$.
\item $\supp b\cap \supp c=\supp d$.
\end{itemize}
\end{Lemma}
Consequently there exist at most $(\weight b)^2$ equivalence  
classes of FK sentences $[x_i]_{i=1}^n$ such that
$b\sim [x_i]_{i=1}^{n-1}$ and $c\sim x_n$. See Figure
3 for an example of two such equivalence classes and
their pictorial description. 
\begin{figure}[h]
\label{fig-twoeqFK}
\begin{center}
\includegraphics[scale=.49]{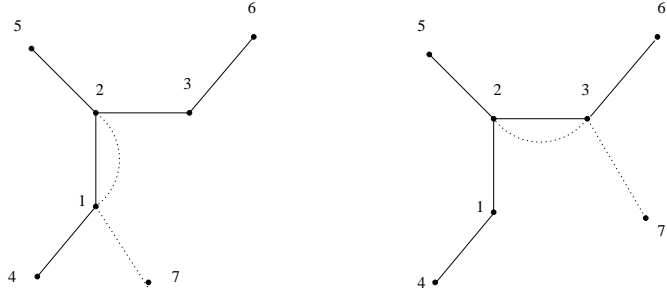}
\caption{Two 
inequivalent FK sentences $[x_1,x_2]$ corresponding
to
$b=141252363$ (solid) and 
$c=1712\sim 3732$ (dashed).}
\end{center}
\end{figure}
\proof Sufficiency is easy to check. We omit the details.
We turn to the proof of necessity. To begin with, since $G_a$ is a tree,
$d$ is the unique
geodesic in $G_c\subset G_a$ joining $\gamma_1$ to $\gamma_\ell$,
and hence is also the unique geodesic in $G_b\subset G_a$ joining
$\gamma_1$
to
$\gamma_\ell$. Now $d$ only visits edges of $G_b$ already
visited by the constituent words of $b$. Therefore we have $E_d\subset
E^1_b$, i.~e.,
$d$ is a walk in $G^1_b$. By
Proposition~\ref{Proposition:Acronym} we have
$E^1_c=E_{\gamma_1\cdots \gamma_r}$. By definition
of an FK sentence we have
$E_b\cap E_c\subset E^1_b\cap E^1_c$. It follows that 
$E_b\cap E_c=E_d$.
Finally, we have
$$\card V_a=1+\card E_a=1+\card E_b +1+\card E_c-1-\card E_d=
\card V_b+\card V_c-\card V_d,$$
and hence, since
$\card V_b+\card V_c-\card V_b\cap V_c=\card V_a$,
 the inclusion 
$V_d\subset V_b\cap V_c$ is in fact an equality.
\qed

\begin{Lemma}\label{Lemma:GammaBound}
Let $\Gamma(k,\ell,m)$ denote the set of equivalence classes of FK
sentences $a=[w_i]_{i=1}^m$ consisting of $m$ words
such that $\sum_{i=1}^m \ell(w_i)=\ell$ and $\weight a=k$.
We have
$$\card \Gamma(k,\ell,m)\leq
2^{\ell-m}\left(\begin{array}{c}
\ell-1\\
m-1
\end{array}\right)k^{2(m-1)}.
$$
\end{Lemma}
\proof 
There are exactly $\left(\begin{subarray}{c}
\ell-1\\
m-1
\end{subarray}\right)$ $m$-tuples of positive integers summing to $\ell$
and hence
by Lemma~\ref{Lemma:FKcounting}
there are at most
$2^{\ell-m}\left(\begin{subarray}{c}
\ell-1\\
m-1
\end{subarray}\right)$  ways to prescribe equivalence
classes of FK words $w_1,\dots,w_m$ subject to the
constraint $\sum_{i=1}^m \ell(w_i)=\ell$.
Now fix FK words $w_1,\dots,w_m$ such that $\sum_{i=1}^m \ell(w_i)=\ell$.
By Lemma~\ref{Lemma:Coding} there exist at most $k^{2(m-1)}$ equivalence
classes of FK sentences $b=[x_i]_{i=1}^m$ such $k=\weight b$
and $w_i\sim x_i$ for $i=1,\dots,m$. The result follows.
\qed

\begin{Lemma} For any FK sentence $a=[w_i]_{i=1}^m$
consisting of $m$ words we have
\begin{equation}\label{equation:EulerFK}
m=\card E_a^1-2\weight a+2+\sum_{i=1}^m \ell(w_i).
\end{equation}
\end{Lemma}
\proof Put $M:=\sum_{i=1}^m\ell(w_i)$. Consider the word
$[\alpha_i]_{i=1}^M$ obtained by concatenating the words of
the sentence
$a$. Consider the list $A=[\{\alpha_i,\alpha_{i+1}\}]_{i=1}^{M-1}$
of unordered pairs of letters. Among the entries of $A$ we find
$2\card E_a-\card
E_a^{1}$ of them that are edges of
$G_a$, while the rest correspond to the $m-1$ ``commas'' in the
sentence
$a$; and moreover, since $G_a$ is a tree, we have $\card E_a=\weight a-1$.
The result follows.  \qed

\begin{Proposition}\label{Proposition:FK}
For all positive integers $n,k$  satisfying 
$n\geq 2k-2$
there are at most
\begin{equation}
\label{eq:290304b}
N_{\rm FK}(n,k):={2^nn^{3(n-2k+2)}}
\end{equation}
equivalence classes of weak Wigner words $w$ such that
$\ell(w)=n+1$ and $\weight w=k$. 
\end{Proposition}
\noindent This is a crude but easy-to-apply
version of the estimate one obtains by exploiting the 
idea of ``coding'' introduced 
by F\"{u}redi and Koml\'{o}s
in \cite{furedi}.
\proof Let $w$ be a weak Wigner word. Let $w'$ be the FK
syllabification of $w$. Let $m$ be the number of words in the sentence
$w'$.  We must have $E^1_{w'}=\emptyset$ lest
there exist an edge of $G_w$ visited only once by $w$
and so we must have $m=\ell(w)-2\weight w+2$ by the preceding lemma. 
Therefore $\card\Gamma(k,n+1,n-2k+3)$
bounds the quantity we wish to estimate, whence
the desired result by Lemma~\ref{Lemma:GammaBound},
after a short further calculation which we omit.
\qed
\subsection{Companion estimate}
To exploit the preceding proposition we need
also to bound $E|\xi(w)|$ for all weak Wigner words $w$
such that $\ell(w)=n+1$ and $k=\weight w$. Fix such a word $w$ now. We
claim that
\begin{equation}\label{eq:300304e}
E|\xi(w)|\leq C(3(n+2-2k))\cdot C(2)^{n/2}\,, \;\mbox{\rm with} \;
C(q):=1\vee\sup_{\alpha,\beta}\max_{m=1}^qE|\xi_{\{\alpha,\beta\}}|^m.
\end{equation}
Consider the graph $G_w=(V_w,E_w)$ and let $\ell$ be the number of
edges of $E_w$ visited exactly twice by $w$. We have by
(\ref{equation:FundamentalBound})  and the H\"{o}lder inequality that
$$E|\xi(w)|\leq C(n-2\ell)C(2)^{\ell}.$$
We have $\card E_w\geq \card V_w-1=k-1$ since $G$ is connected,
$n\geq 3\cdot(\card E_w-\ell)+2\ell$ by counting, and hence
$$n-2\ell\leq 3(n+2-2k).$$
The desired estimate now follows since $C(q)$ is a nondecreasing
function of $q$ bounded below by $1$.

\section{Bracelets, polarizations and enumeration}
\label{sec:bracelets}
We have already seen in Section \ref{section:Limcal} that limiting
variances
are determined by the enumeration of CLT word-pairs. In the  current
section, 
we study the structure of such word-pairs and their associated
graphs. These turn out to be classified
by certain ``bracelets with pendant trees''.

\subsection{Graph-theoretical definitions}
\subsubsection{Bracelets}
We say that a graph
$G=(V,E)$ is a {\em bracelet} if there exists an enumeration
$\alpha_1,\dots,\alpha_r$ of $V$ such that
$$E=\left\{\begin{array}{rl}
\{\{\alpha_1,\alpha_1\}\}&\mbox{if $r=1$,}\\
\{\{\alpha_1,\alpha_2\}\}&\mbox{if $r=2$,}\\
\{\{\alpha_1,\alpha_2\},\{\alpha_2,\alpha_3\},
\{\alpha_3,\alpha_1\}\}&\mbox{if
$r=3$,}\\
\{\{\alpha_1,\alpha_2\},\{\alpha_2,\alpha_3\},
\{\alpha_3,\alpha_4\},\{\alpha_4,\alpha_1\}\}&\mbox{if
$r=4$,}\\
\end{array}\right.
$$
and so on.
We call $r$ the {\em circuit length} of the
bracelet $G$. 

\subsubsection{Unicyclic graphs}
We say that a graph $G=(V,E)$ is {\em unicyclic}
if $G$ is connected and $\card V=\card E$. In other words, a unicyclic
graph is a connected graph with one too many edges to be a tree. Any
bracelet  of circuit length 
$\neq 2$ is
unicyclic. However, a bracelet of circuit length $2$ is a
tree.

\begin{Proposition}\label{Proposition:FindTheBracelet} Let
$G=(V,E)$ be a unicyclic graph. For each edge
$e\in E$ put
$G\setminus e=(V,E\setminus\{e\})$.
 Let $Z$ be the subgraph of $G$ consisting of all $e\in E$
such that $G\setminus e$ is connected,  along with all attached vertices.
Let $r$ be the number of edges of $Z$.
 Let
$F$ be the graph obtained from $G$ by  deleting all edges of $Z$. The
following statements hold:
\begin{enumerate}
\item $F$ is a forest with exactly $r$ connected components.
\item If $G$ has a degenerate edge, then $r=1$.
\item If $G$ has no degenerate edge, then $r\geq 3$.
\item $Z$ meets each connected component of $F$ in exactly one vertex.
\item $Z$ is a bracelet  of circuit length $r$.
\item For all $e\in E$ the following conditions are equivalent:
\begin{enumerate}
\item $G\setminus e$ is connected.
\item $G\setminus e$ is a tree.
\item $G\setminus e$ is a forest.
\end{enumerate}
\end{enumerate}
\end{Proposition}
We call $Z$ the {\em bracelet} of $G$. We call $r$ 
the {\em circuit length}
of $G$, and each of the components of $F$ we call a {\em pendant tree}. 
\proof The proposition is well-known in principle. We just explain how
to prove statement 5 and omit the remaining details. Pick an edge
$e=\{\alpha,\beta\}$ of $G$ so that $G\setminus e$ is a spanning tree.
Then $e$ is an edge of $Z$, and it is not difficult to verify that the
edges
of
$Z$ distinct from $e$ are the edges of the tree $G\setminus e$ visited by
the unique geodesic in
$G\setminus e $ joining $\alpha$ to $\beta$. So it is clear that $Z$
is a bracelet.
\qed
\subsection{The bracelet of a CLT
word-pair}
\label{subsection:CLTWordPairBracelets}
Fix a CLT word-pair $[w,x]$. Let $G=G_{[w,x]}$ be the
associated graph.
\subsubsection{} By
Proposition~\ref{Proposition:CLTWordPairMomentControl}, either $G$ is
unicyclic or $G$ is a tree. If $G$ is unicyclic, we define the {\em
bracelet}, {\em circuit length}, and 
{\em pendant trees} of $[w,x]$ to be the same as those
defined for
$G$ by Proposition~\ref{Proposition:FindTheBracelet}.
Suppose now that $G$ is a tree. Then there exists by
Proposition~\ref{Proposition:CLTWordPairMomentControl}(3)(c) a unique edge
of
$G$ visited exactly twice by $w$ and twice by $x$; 
this edge and attached
vertices we declare to be the {\em bracelet} of $[w,x]$, and we declare
the {\em circuit length} of
$[w,x]$ to be the circuit length of its bracelet, namely $2$.
Erasing the edge of the bracelet from $G$ leaves a forest of two
components; as before, we call the components {\em
pendant trees} of $[w,x]$.
Note that in all cases the circuit length of $[w,x]$ 
depends only on the equivalence
class of the word-pair $[w,x]$. 
\subsubsection{}
Let $Z$ and $r$ be the bracelet and circuit length of $[w,x]$,
respectively. Note that 
$G$ is unicyclic or a tree according to whether $r\neq 2$ or $r=2$.
\begin{figure}[h]
\begin{center}
\includegraphics[scale=.49]{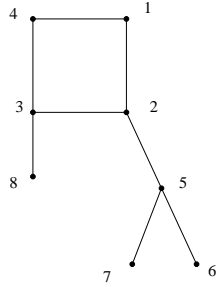}
\caption{The bracelet 1234 of circuit length 4,
and the pendant trees, associated with the CLT word-pair
$[12565752341,2383412]$}
\end{center}
\end{figure}

\subsubsection{}\label{sec-almost2}
Now write $w=[\alpha_i]_{i=1}^{\ell(w)}$ and
$x=[\beta_j]_{j=1}^{\ell(x)}$. Let $\check{w}$ and $\check{x}$ be the
words obtained by dropping the last letters of $w$ and $x$, respectively.
Let $\sigma$ be any cyclic permutation of $\{1,\dots,\ell(\check{w})\}$
and let $\tau$ be any cyclic permutation of $\{1,\dots,\ell(\check{x})\}$.
Then $[\check{w}^\sigma
\alpha_{\sigma(1)},
\check{x}^\tau\beta_{\tau(1)}]$ is again a CLT word-pair with
associated graph, bracelet and circuit length the same as for
$[w,x]$. (The ``exponential notation'' used here was defined in
\S\ref{subsubsection:CyclicPermutations}.) 
We declare the ordered pair
$(\sigma,\tau)$ to be a {\em polarization} of $[w,x]$ if the last edge
of $G$ visited by the walk 
$\check{w}^\sigma
\alpha_{\sigma(1)}$ equals the last edge of $G$ visited by the
walk $\check{x}^\tau\beta_{\tau(1)}$. Note that the set of polarizations of
a CLT word-pair depends only on its equivalence class. 
The notions of bracelet and polarization are linked by the
following result.

\begin{Lemma}\label{Lemma:EasyBraceletCharacterization}
Let $[w,x]$ be a CLT word-pair. Put $G=G_{[w,x]}$.
Let $Z$ and $r$ denote the bracelet and circuit length of $[w,x]$,
respectively. Let
$e$ be an edge of
$G$. Then: (i)
$e$ is an edge of
$Z$ if and only if both words/walks $w$ and $x$ visit $e$. (ii)
Unless $r=2$, there exist exactly $r$ polarizations of $[w,x]$; but if
$r=2$, there exist exactly $4$ polarizations of $[w,x]$.
\end{Lemma}
In the example of Figure 4, the four polarizations lead to the
CLT word-pairs 
$$\begin{array}{ll}
\left[12565752341,1238341\right]\,,& \left[25657523412,2383412
\right]\,,\\
\left[34125657523,3834123
\right]\,,& \left[41256575234,4123834\right]\,.
\end{array}$$
\proof
If $G$ is a tree, part (i) of
the lemma holds by definition of $Z$, while part (ii)
is a consequence of 
Proposition \ref{Proposition:CLTWordPairMomentControl}(3).
So assume for the rest of the proof that $G$ is unicyclic. 
By Proposition~\ref{Proposition:CLTWordPairMomentControl}(4)(b)
each edge of $G$ is visited a total of exactly 
two times by $w$ and $x$,
and so part (ii) of the proposition follows immediately from part (i).
We have only to prove  part (i).
($\Rightarrow$) The graph
$G\setminus e$ obtained by deleting $e$ is by hypothesis  a tree. If one
of the words/walks $w$ or $x$ fails to visit $e$, say the former, then by
the parity principle 
 $w$ must visit every edge of
$G$ an even number of times. But then, due
to Proposition \ref{Proposition:CLTWordPairMomentControl}(4)(b),
it is impossible for $w$ to visit
any edge of
$G$ visited by
$x$, which is a contradiction.  ($\Leftarrow$)  
By hypothesis and 
Proposition \ref{Proposition:CLTWordPairMomentControl}(4)(b)
the walk
$w$ visits $e$ exactly once,  hence some cyclic
permutation of
$\check{w}$ is a walk on
$G\setminus e$ the set of endpoints of which equals $e$, hence
$G\setminus e$ is connected, and hence $e$ is an edge of $Z$. 
\qed

\begin{Lemma}\label{Lemma:Zip}
Let $G$ be a forest. Let $w$ and $x$ be words both
admitting interpretation as walks on $G$.
Assume that jointly $w$ and $x$ visit every
edge of
$G$ either exactly twice or never. 
(Necessarily then, both $w$
and
$x$ are FK words.) 
Assume further that
$w$ and
$x$ have at least one letter in common. Then exactly one of the following
conditions holds:
\begin{enumerate}
\item  $w$ and $x$ have acronyms which are either equal or mirror images,
but have no letters in common apart from those shared by their
acronyms.
\item $w$ and $x$ are Wigner words with exactly
one letter in common, but this common letter does not appear
as the first letter of both words.
\end{enumerate}
\end{Lemma}
\proof 
Since any subgraph of a forest is again a forest, we may
assume without loss of generality that
$$G=(V,E)=(V_w\cup V_x,E_w\cup E_x).$$
Since $w$ and $x$ have at least one letter in common, in fact $G$ is a
tree. Let
$w_*$ and
$x_*$ be the acronyms of $w$ and $x$, respectively. Note
that $w_*$ (resp., $x_*$) is the unique geodesic in $G$ with the same
initial and terminal vertices as $w$ (resp., $x$). By the parity
principle and the hypotheses we have
$$\begin{array}{cl}
&\{e\in E\vert \mbox{$w$ visits $e$ exactly once}\}\\
=&\{e\in E\vert \mbox{$w_*$
visits
$e$}\}
=E_w\cap E_x
=\{e\in E\vert
\mbox{$x_*$ visits $e$}\}\\
=&\{e\in E\vert \mbox{$x$ visits $e$ exactly once}\},
\end{array}
$$
hence $w_*$ and $x_*$ are words of the same length,
say $\ell$, and we have
$$\ell=1+\card
E_w\cap E_x.$$ If $\ell>1$, then the words $w_*$ and $x_*$ must 
either be equal or mirror images of each other. 
If $\ell=1$, then $w$ and $x$ are Wigner words since each visits every
edge of $G$ either exactly twice or never, but note that we need not in
this case have equality of $w_*$ and $x_*$. Finally,
since
$G$,
$G_w$ and
$G_x$ are trees, we have
$$\card V=1+\card E=
1+\card E_w+\card E_x-\card E_w\cap E_x
=\card V_w+\card V_x-\ell,$$
which finishes the proof. 
\qed

\begin{Proposition}\label{Proposition:Spin}
Fix closed words $w$ and $x$ each of length $\geq 2$.
Put $k=\ell(\check{w})$ and $\ell=\ell(\check{x})$.
Let $\sigma$ (resp., $\tau$) be a cyclic permutation
of $\{1,\dots,k\}$ (resp., $\{1,\dots,\ell\}$).
The following statements are equivalent:
\begin{enumerate}
\item $[w,x]$ is a CLT word-pair of which 
$(\sigma,\tau)$ is a polarization.
\item $\check{w}^\sigma$ and $\check{x}^\tau$ are FK words with acronyms
either equal or mirror images, and with no letters in common apart from
those shared by their acronyms.
\end{enumerate}
\end{Proposition}
\noindent We remark that under the equivalent conditions above,
the common length of the acronyms of $\check{w}^\sigma$ and
$\check{x}^\tau$
equals the circuit length of $[w,x]$. 

\proof The implication 2$\Rightarrow$1 is easy to check. We omit the
details. We turn directly to the proof of the implication
1$\Rightarrow$2.
Write $w=[\alpha_i]_{i=1}^{k+1}$ and $x=[\beta_j]_{j=1}^{\ell+1}$. Let
$e$ be the last edge of
$G=G_{[w,x]}$ visited by the walks $\check{w}^\sigma \alpha_{\sigma(1)}$ 
and $\check{x}^\tau \beta_{\tau(1)}$. Note that 
$e$ by
Lemma~\ref{Lemma:EasyBraceletCharacterization} is automatically an edge
of the bracelet of $[w,x]$. Unless $r=2$, let
$G'$ be the graph obtained by deleting
$e$ from
$G$, but if $r=2$ put $G'=G$. Then in all cases $G'$ is a tree,
and the words $\check{w}^\sigma$ and $\check{x}^\tau$ are
walks on $G'$ satisfying the hypotheses of
Lemma~\ref{Lemma:Zip}.
Were $\check{w}^\sigma$ and
$\check{x}^\tau$ to be Wigner words with exactly one letter in common not
appearing as the first letter of both words, the graph
$G$ would have two degenerate edges, which by
Proposition~\ref{Proposition:CLTWordPairMomentControl} is impossible.
\qed

\subsection{Enumeration of CLT
word-pairs by Wigner words}\label{subsection:FurtherEnumerations} We are
now ready to state an enumeration formula for CLT word-pairs
similar to the enumeration formulas 
(\ref{equation:WignerAutoEnumeration}) and (\ref{equation:Marriage}),
albeit rather more complicated.

\subsubsection{Enumerative
apparatus}\label{subsubsection:EnumerativeApparatus} Let
$[\gamma_i]_{i=1}^\infty$ be a sequence of distinct letters.
For each positive integer $i$ choose cross sections $U_i$ and $V_i$ of the
set of Wigner words. Make these choices so as to achieve the following
conditions:
\begin{itemize}
\item For all $i$, every word belonging to $U_i\cup V_i$ begins with
$\gamma_i$, but no word belonging to $U_i$ has a letter other than
$\gamma_i$ in common  with any word belonging to $V_i$. 
\item For all distinct $i$ and $j$, every word belonging to
$U_i\cup V_i$ is disjoint from every word belonging to $U_j\cup V_j$.
\end{itemize}
Let $\varphi$ be a real-valued function 
defined for all sentences.
Assume that $\varphi(a)$ depends only on the equivalence class of
$a$ and vanishes when the sum of the lengths of the constituent words of
$a$ is sufficiently large, in which case the support of $\varphi$
consists of only finitely many equivalence classes of sentences.
\subsubsection{Enumeration of CLT word-pairs}
We have
\begin{equation}\label{equation:CLTWordPairEnumeration}
\begin{array}{cl}
&\displaystyle\sum_{a} \varphi(a)\\\\
=&\displaystyle\sum_{r=1}^\infty
\sum_{u_1\in U_1}
\cdots \sum_{u_r\in U_r}
\sum_{v_1\in V_1}
\cdots \sum_{v_r\in V_r}
\sum_\sigma \sum_\tau\\\\
&
\left\{
\begin{array}{rl}
\varphi([u^\sigma \alpha_{\sigma(1)},v^\tau
\beta_{\tau(1)}])&\mbox{if
$r=1$,}\\
\left(\varphi([u^\sigma \alpha_{\sigma(1)},v^\tau
\beta_{\tau(1)}])+\varphi([u^\sigma \alpha_{\sigma(1)},\bar{v}^\tau
\bar{\beta}_{\tau(1)}])\right)/4&\mbox{if $r=2$,}\\
\left(\varphi([u^\sigma \alpha_{\sigma(1)},v^\tau
\beta_{\tau(1)}])+\varphi([u^\sigma \alpha_{\sigma(1)},\bar{v}^\tau
\bar{\beta}_{\tau(1)}])\right)/r&\mbox{if $r\geq 3$,}\\
\end{array}\right.
\end{array}
\end{equation}
where:
\begin{itemize}
\item $a$ ranges over any cross-section of the set of CLT word-pairs;
\item $u=u_1\cdots u_r=[\alpha_i]_{i=1}^{\ell(u)}$;
\item $v=v_1\cdots v_r=[\beta_i]_{i=1}^{\ell(v)}$ and $\bar{v}=v_r\cdots
v_1=[\bar{\beta}_i]_{i=1}^{\ell(v)}$;
\item $\sigma$ ranges over cyclic permutations of $\{1,\dots,\ell(u)\}$;
and
\item $\tau$ ranges over cyclic permutations of
$\{1,\dots,\ell(v)\}$.
\end{itemize}
One verifies that there is neither under- nor over-counting by applying
Proposition~\ref{Proposition:Acronym} (which gives the structure of FK
words) and Proposition~\ref{Proposition:Spin} (which gives the structure
of CLT word-pairs) in a straightforward way. We omit further details.

\section{Proof of Theorem~\ref{Theorem:Main}}
\label{sec:proofII}
\subsection{Further generating functions}
Fix a sentence 
$$a=[w_i]_{i=1}^n=[[\alpha_{ij}]_{j=1}^{\ell(w_i)}]_{i=1}^n$$
consisting of $n$ words.
\subsubsection{}
 Let $t=[t_i]_{i=1}^n$
be an
$n$-tuple of independent (algebraic) variables and put
$$H(a,t)=\sum_p
 H_{p}(a)\prod_{i=1}^n t_i^{p_i+\ell(w_i)}
$$
where $p=[p_i]_{i=1}^n$ ranges over $n$-tuples of (nonnegative) integers.
We view $H(a,t)$ as a formal power
series in $t_1,\dots,t_n$ with random variable coefficients, not as an
analytic function of $t$. In other words,
$H(a,t)$ is just a device for manipulating the
infinite array $[H_p(a)]$ of random variables.
We  write
$$MH(a,t)=M(a)H(a,t),\;\;\;\overline{M}H(a,t)=\overline{M}(a)H(a,t)$$
in order to abbreviate notation.

\subsubsection{} Unraveling the definition
of $H(\cdot,\cdot)$ in the case of a single word
$w=[\alpha_j]_{j=1}^{\ell(w)}$, we find that
\begin{equation}
\label{friday1}
H(w,t)=t^{\ell(w)} 
\sum_{\pi}
\prod_{j=1}^{\ell(w)} [D(\kappa(\alpha_j))t]^{\pi_j}=
\prod_{j=1}^{\ell(w)}\frac{t}{1-tD(\kappa(\alpha_j))}\,,
\end{equation} 
where $\pi=[\pi_j]_{j=1}^{\ell(w)}$ ranges over $\ell(w)$-tuples
of nonnegative integers.  From
(\ref{friday1}), it follows that
\begin{eqnarray*}
H(w\alpha_1,t)&=&
\left(\frac{t}{1-tD(\kappa(\alpha_1))}\right)^2\cdot
\prod_{j=2}^{\ell(w)}
\frac{t}{1-tD(\kappa(\alpha_j))}\\
&=&\left(t^2\frac{d}{dt}\frac{t}{1-tD(\kappa(\alpha_1))}\right)\cdot
\prod_{j=2}^{\ell(w)}
\frac{t}{1-tD(\kappa(\alpha_j))}.
\end{eqnarray*}
Taking the sum over all
cyclic permutations $\sigma$ of $\{1,\ldots,\ell(w)\}$,
and arguing similarly, we find that
\begin{equation}
\label{friday2}
\sum_{\sigma}
H(w^\sigma\alpha_{\sigma(1)},t)=t^2 \frac{\partial}{\partial t}
H(w,t)\,.
\end{equation}

\subsubsection{}
Returning now to the general situation, from (\ref{friday1}) we get 
the identity
\begin{equation}\label{equation:HGeneration}
H(a,t)=\prod_{i=1}^n\prod_{j=1}^{\ell(w_i)}
\frac{t_i}{1-D(\kappa(\alpha_{ij}))t_i}=\prod_{i=1}^n H(w_i,t_i)\,.
\end{equation} From
(\ref{friday2}) and (\ref{equation:HGeneration})
we get the differentiation formula
\begin{equation}\label{equation:HDifferentiationFormula}
t_1^2\frac{\partial}{\partial t_1}
\cdots t_n^2\frac{\partial}{\partial t_n}H(a,t)=
\sum_{\sigma_1}\cdots \sum_{\sigma_n}
H([w_i^{\sigma_i} \alpha_{i,\sigma_i(1)}]_{i=1}^n,t),\end{equation}
where in the sum $\sigma_i$ ranges over
cyclic permutations of $\{1,\dots,\ell(a_i)\}$.
We emphasize that these identities are to be interpreted formally, i.~e.,
all the expressions are to be expanded as power series in
$t_1,\dots,t_n$ in evident fashion and then coefficients of like
monomials in the
$t_i$ are  to be equated.

\subsubsection{}
For each Wigner word $w$ we define
$$\Phi^{(w)}(c,t)=\sum_{p=0}^\infty \Phi^{(w,p)}(c)t^{\ell(w)+p}.$$
As with the generating functions introduced above, this, too, is to be
viewed as formal power series in $t$.
By Lemma~\ref{Lemma:ConditionalPhiInterpretation}
we have
\begin{equation}\label{equation:ConditionalPhiInterpretation1}
E( MH(w,t)
\vert\kappa(\alpha))=\Phi^{(w)}(\kappa(\alpha),t)\;\;\;\mathrm{a.s.}
\end{equation}
where to make sense of formula, both sides are expanded in
powers of $t$, the integrals on the left are computed term by term,
and then coefficients of like powers of $t$ are to be set equal a.s.
By Lemma~\ref{Lemma:PhiDecomposition}
we have
\begin{equation}\label{equation:PhiDecomposition1}
\Phi(c,t)=\sum_w \Phi^{(w)}(c,t)
\end{equation}
where $w$ ranges over a cross-section of the
set of Wigner words.  Note that in the sum on the right, for every fixed
degree
$n$, there are only finitely many terms in which the coefficient of $t^n$
is nonvanishing.

\begin{Lemma}\label{equation:LemmaKeyCovariance}
We have an identity
\begin{equation}\label{equation:CovarianceFormulaBis}
\sum_{i=1}^\infty \sum_{j=1}^\infty
EY_iY_j\cdot x^{i}y^{j}=
x\frac{\partial}{\partial x} y\frac{\partial}{\partial y}
\bigg(
2\Theta(x,y)+\Psi(x,y)\bigg)
\end{equation}
of formal power series, where $[Y_i]_{i=1}^\infty$
is the Gaussian family defined in Lemma~\ref{Lemma:GaussMarriage}.
\end{Lemma}
\proof Let 
$[\gamma_i,U_i,V_i]_{i=1}^\infty$
be the enumerative apparatus introduced in
\S\ref{subsection:FurtherEnumerations}. In anticipation of applying
enumeration formula (\ref{equation:CLTWordPairEnumeration}) we temporarily
``freeze'' data specifying a single term on the right side of that
formula:
\begin{itemize}
\item Let $r$ be a positive integer.
\item Let $u_1\in U_{1},\dots,u_r\in U_{r}$
and $v_1\in V_{1},\dots,v_r\in V_{r}$.
\item Let $u=[\alpha_i]_{i=1}^{\ell(u)}=u_1\cdots u_r$.
\item Let
$v=[\beta_j]_{j=1}^{\ell(v)}$ be equal either to
$v_1\cdots v_r$ or to  $v_r\cdots v_1$.
\item Let $\sigma$ be a cyclic permutation of $\{1,\dots,\ell(u)\}$.
\item Let $\tau$ be a
cyclic permutation of $\{1,\dots,\ell(v)\}$.
\end{itemize}
By direct appeal to the definitions we have
\begin{equation}
\label{equation:OverlineMFactorization}
\begin{array}{cl}
&\overline{M}([u\alpha_1,v\beta_1])=\overline{M}([u^\sigma
\alpha_{\sigma(1)},v^\tau \beta_{\tau(1)}])\\\\
=&\displaystyle\prod_{i=1}^r(M(u_i)M(v_i))\cdot
\left\{\begin{array}{ll}
d^{(2)}(\kappa(\gamma_1))&
\mbox{if $r=1$,}\\
s^{(4)}(\kappa(\gamma_1),\kappa(\gamma_2))-
s^{(2)}(\kappa(\gamma_1),\kappa(\gamma_2))^2&\mbox{if
$r=2$,}\\
K_r(\kappa(\gamma_1),\dots,\kappa(\gamma_r))&\mbox{if $r\geq 3$.}
\end{array}\right.
\end{array}
\end{equation}
To understand this formula, notice that the right side is a product 
of factors associated to pendant trees times a factor arising from the
bracelet. Put 
$$\FFF=\sigma([\kappa(\gamma_i)]_{i=1}^\infty).$$ 
We then have the following identities:
\begin{equation}
\begin{array}{cl}
&\displaystyle
\sum_\sigma\sum_\tau
E(\overline{M}H([u^\sigma \alpha_{\sigma(1)},
v^\tau \beta_{\tau(1)}],[x,y])\vert\FFF)\\\\
=&\displaystyle
\sum_\sigma\sum_\tau
E(\overline{M}([u\alpha_1,v\beta_1])H(u^\sigma \alpha_{\sigma(1)},x)
H(v^\tau \beta_{\tau(1)},y)\vert\FFF)\\\\
=&\displaystyle
x^2\frac{\partial}{\partial x} y^2\frac{\partial}{\partial y}
E(\overline{M}([u\alpha_1,v\beta_1])H(u,x)H(v,y)\vert\FFF)\\\\
=&\displaystyle  x^2\frac{\partial}{\partial x} 
\prod_{i=1}^r
\Phi^{(u_i)}(\kappa(\gamma_i),x)
\cdot y^2\frac{\partial}{\partial
y}\prod_{i=1}^r\Phi^{(v_i)}(\kappa(\gamma_i),y)\\\\ &\cdot
\left\{\begin{array}{ll}
d^{(2)}(\kappa(\gamma_1))&
\mbox{if $r=1$,}\\
s^{(4)}(\kappa(\gamma_1),\kappa(\gamma_2))-
s^{(2)}(\kappa(\gamma_1),\kappa(\gamma_2))^2&\mbox{if
$r=2$,}\\
K_r(\kappa(\gamma_1),\dots,\kappa(\gamma_r))&\mbox{if $r\geq 3$.}
\end{array}\right.
\end{array}
\end{equation}
Here all the conditional expectations are to be calculated
by expanding formally in powers of $x$ and $y$ and then integrating term by
term; in the same spirit the equal signs are to be interpreted as a.s.\
equality term by term between formal power series. The
preceding holds at  the second equality by the differentiation formula
(\ref{equation:HDifferentiationFormula}), and at the third equality by
(\ref{equation:ConditionalPhiInterpretation1}). Now take expectations
(again, integrating term by term), and then apply identity
(\ref{equation:PhiDecomposition1}) and enumeration formula
(\ref{equation:CLTWordPairEnumeration}) to find that
\begin{equation}\label{equation:CovarianceFormulaTer}
\sum_{[u,v]}\overline{M}H([u,v],[x,y])=
x^2\frac{\partial}{\partial x} y^2\frac{\partial}{\partial y}
\bigg(
2\Theta(x,y)+\Psi(x,y)\bigg)
\end{equation}
where on the left $[u\;v]$ ranges over a cross-section of the set of CLT
word-pairs. The result now follows by definition
of the random variables $Y_i$. 
\qed

\subsection{End of the proof of Theorem~\ref{Theorem:Main}}
By limit formula (\ref{equation:SentenceExpansion3}),
 Lemma~\ref{Lemma:GaussMarriage} and 
Lemma~\ref{equation:LemmaKeyCovariance}, 
we have for every nonnegative integer $n$ that
$$\lim_{k\rightarrow\infty} EZ_{f,k}^n=EY_f^n$$
where 
$$Y_f=\sum_{i=1}^\infty \langle t^i,f(t)\rangle Y_i.$$
So the method of moments gives the result.
\qed

\section{Proof of Theorem~\ref{Theorem:MainBis}}\label{sec:CalcMean}
Because of the strong similarity between the proofs of
Theorem~\ref{Theorem:Main} and Theorem~\ref{Theorem:MainBis},
and because of the tedious nature of latter proof (the
necessary enumerations are rather involved),  we
proceed quite rapidly, omitting many details. But we strive to provide
all the important ``landmarks'' so that the reader won't get lost.

\subsection{Further random variables indexed by sentences}
\label{subsection:SentenceIndexedBis} 
We enlarge the supply of random
variables introduced in \S\ref{subsection:SentenceIndexed},
as follows. 
\subsubsection{}
Given any finite nonempty set $\NNN$ of letters, let
$\{\kappa_\NNN(\alpha)\}$ be a letter-indexed color-valued family of
i.i.d.\ random variables with common law $\theta_\NNN$. Then, with
$\NNN$ as above, for any word $w$ and integer $p$ we define
$M_\NNN(w)$ and $H_{p,\NNN}(w)$  by repeating
the definitions of $M(w)$ and $H_p(w)$,  see
\S\ref{sec-almost0} and \S\ref{sec-almost1},
with $\theta_\NNN$
in place of $\theta$. In fact, these random variables
were already considered in the course of the proof of
Theorem~\ref{Theorem:Preliminary}, 
see equation (\ref{equation:WordExpansion5}).

\subsubsection{}
Given distinct letters $\alpha$ and $\beta$,
let  $[\beta\mapsto \alpha]$ be the unique map
of letter space to itself sending $\beta$ to $\alpha$
but fixing all other letters. 
Given also a word $w$, let $[\beta\mapsto \alpha]_*w$ be the word
obtained by applying $[\beta\mapsto \alpha]$ letter by letter
to
$w$.

\subsubsection{}
Let $\psi$ be any map of letter space to itself.
Let $w=[\alpha_i]_{i=1}^n$ be any closed word.
Put
$$
M(w,\psi)
=\prod_{\begin{subarray}{c}
e=\{\alpha,\beta\},\\
\mbox{\tiny
edge of
$G_w$}
\end{subarray}}
\left\{\begin{array}{cl}
0&\mbox{if $\nu(e)=1$}\\
s^{(\nu(e))}(\kappa(\psi(\alpha)),
\kappa(\psi(\beta)))&\mbox{if $\nu(e)>1$ and $\alpha\neq \beta$,}\\
d^{(\nu(e))}(\kappa(\psi(\alpha)))
&\mbox{if $\nu(e)>1$ and $\alpha=\beta$,}
\end{array}\right.
$$ 
where $\nu(e)$ is the number of visits made by $w$ to $e$.
Note that if $\psi$ is the identity map, then $M(w,\psi)=M(w)$.
The only case of the generalization
$M(\cdot,\cdot)$ of $M(\cdot)$ figuring in our limit formulas is that
in which
$w$ is a Wigner word and
$\psi=[\beta\mapsto
\alpha]$ for some distinct letters $\alpha$ and $\beta$ appearing in
$w$. Note that in that case $M(w,\psi)$ depends only on $s^{(2)}$,
not on
$\{s^{(k)}\}_{k\neq 2}\cup \{d^{(k)}\}$.

\subsection{Approximation of $\langle
\overline{L}(\NNN),x^n\rangle$ at CLT scale}
Fix a positive integer $n$. 
Let $[\NNN_k]_{k=1}^\infty$ be as in Assumption \ref{ass-main}.
Starting again with formula (\ref{equation:WordExpansion2}),
it is possible to
obtain the formula
\begin{equation}\label{equation:CLTScaleBis}
\begin{array}{cl}
&\displaystyle\lim_{k\rightarrow\infty}
 N_k \cdot \left(
\langle
\overline{L}(\NNN_k),x^n\rangle-
 \sum_w
EM_{\NNN_k}H_{n+1-\ell(w),\NNN_k}
(w)\right)\\\\
=&\displaystyle
-\frac{1}{2}\sum_{[u,\alpha,\beta]}EM(u,[\beta\mapsto
\alpha])H_{n+1-\ell(u)}([\beta\mapsto \alpha]_*u)\\\\
&\displaystyle +\sum_v
EMH_{n+1-\ell(v)}(v)
\end{array}
\end{equation}
where:
\begin{itemize}
\item $w$ ranges over a cross-section of the set of Wigner words;
\item $[u,\alpha,\beta]$ ranges over a cross-section of the set of marked
Wigner words; and
\item $v$ ranges over a cross-section of the set of critical weak Wigner
words.
\end{itemize}
Note that only finitely many nonzero terms appear in the sums.
Since the proof of (\ref{equation:CLTScaleBis})
is quite similar to that of
(\ref{equation:WordExpansion3}), if rather more complicated, we omit the
details. We only give the following hint to the reader. Let us
return for a moment to the set up of
\S\ref{subsubsection:ApproximateIntegrationExample}.
We have
$$N(S(\NNN,w)-EM_{\NNN}H_{n+1-\ell(w),\NNN}
(w))=-N^{1-r}\sum_{\begin{subarray}{c}
(\beta_1,\dots,\beta_r)\in\NNN^r\\
\card\{\beta_1,\dots,\beta_r\}<r
\end{subarray}}f(\kappa_0(\beta_1),\dots,\kappa_0(\beta_r))$$
and up to an $O(N^{-1})$ error the right side equals
$$-N^{-r}\sum_{1\leq i<j\leq r}
\sum_{
(\beta_1,\dots,\beta_r)\in\NNN^r}
f(\kappa_0(\beta_{[j\mapsto i](1)}),\dots,
\kappa_0(\beta_{[j\mapsto i](r)}))$$
where $[j\mapsto i]$ denotes the map of $\{1,\dots,r\}$
to itself sending $j$ to $i$ and fixing all other elements. 
In the case that color space consists of a single color, 
the preceding remark boils down to the observation that
$$N(N-1)\cdots (N-r+1)-N^{r}=-\left(\begin{array}{c}
r\\
2
\end{array}\right)N^{r-1}+\cdots$$
where the omitted terms are $O(N^{r-2})$.
\subsection{Enumeration of marked Wigner words by Wigner words}
We use the enumerative apparatus introduced in 
\S\ref{subsubsection:EnumerativeApparatus}. We have
\begin{equation}\label{equation:MarkedWignerEnumeration}
\sum_{[w,\alpha,\beta]}
\varphi([w,\alpha,\beta])=\sum_{r=1}^\infty
\sum_{u_1\in U_1}\cdots \sum_{u_{r}\in U_{r}}
\sum_{v_1\in V_1}\cdots \sum_{v_r\in V_r}\sum_\sigma
\varphi([u^\sigma \alpha_{\sigma(1)},\gamma_1,\gamma_{r+1}])
\end{equation}
where:
\begin{itemize}
\item $[w,\alpha,\beta]$ ranges over any cross-section of the set of
marked Wigner words;
\item $u=u_1\cdots u_{r}([\gamma_1\mapsto \gamma_{r+1}]_*v_1)v_r\cdots
v_2=[\alpha_i]_{i=1}^{\ell(u)}$ and
\item $\sigma$ ranges over cyclic permutations of $\{1,\dots,\ell(u)\}$.
\end{itemize}
Note that in this setting
\begin{equation}\label{equation:MarkedWignerEnumerationM}
M(u^\sigma \alpha_{\sigma(1)},[\gamma_{r+1}\mapsto \gamma_1])=
\prod_{i=1}^r (M(u_i)M(v_i))\cdot K_r(\kappa(\gamma_1),\dots,
\kappa(\gamma_r)).
\end{equation}
The intuition behind (\ref{equation:MarkedWignerEnumeration})
is as follows. Let $[w,\alpha,\beta]$ be a marked Wigner word, write
$w=[\alpha_i]_{i=1}^{\ell(w)}$, and let $\check{w}$
be the result of dropping the last letter of $w$.
After replacing $w$ by $\check{w}^\sigma \alpha_{\sigma(1)}$
for a certain uniquely determined cyclic permutation
$\sigma$ of
$\{1,\dots,\ell(\check{w})\}$, we may assume that
$\alpha$ is the first letter of $w$ and that every appearance of
$\alpha$ in
$\check{w}$ precedes every appearance of
$\beta$. We may then view
$w$ as a walk out and back on the geodesic connecting $\alpha$ to
$\beta$ in the tree $G_w$ punctuated by sidetrips on the trees
hanging from that geodesic. More precisely, an argument employing 
Proposition~\ref{Proposition:ExtremalWeakWigner} (which gives the
structure of Wigner words),
Proposition~\ref{Proposition:Acronym} (which gives the structure
of FK words) and Lemma~\ref{Lemma:Zip} shows that
there is neither under- nor over-counting 
in (\ref{equation:MarkedWignerEnumeration}).
We omit the details.

\subsection{The bracelet of a critical weak Wigner
word} Let $w=[\alpha_i]_{i=1}^{\ell(w)}$ be a critical weak
Wigner word. Put $G=(V,E)=G_w=(V_w,E_w)$.
\label{subsection:BraceletDef}
\subsubsection{}  According to
Proposition~\ref{Proposition:CriticalMomentOkay},
either $G$ is unicyclic or $G$ is a tree. If $G$ is unicyclic,
then we define the {\em bracelet} and {\em circuit length} of $w$
to be the same as defined for $G$ in
Proposition~\ref{Proposition:FindTheBracelet}. 
If $G$ is a tree, then there exists a unique edge $e$ of $G$
visited exactly $4$ times by $w$; this edge and attached vertices
we declare to be the {\em bracelet} of $w$, and we declare the {\em
circuit length} of $w$ to be that of its bracelet, namely $2$.  

\subsubsection{}
Let $Z$ and $r$ be the bracelet and circuit length of $w$, respectively.
Note that $r\neq 2$ or $r=2$ according to whether $G$ is 
unicyclic or 
a tree.
Note that in all cases the graph obtained from $G$ by deleting the edges
of $Z$ is a forest with exactly $r$ connected components each of which
meets the bracelet in exactly one vertex; again we have a
picture of ``bracelet with pendant trees''. Note that in all cases
$w$ makes a total of
$2r$ visits to edges of
$Z$. Moreover the walk $w$ visits each edge of the
bracelet exactly twice, unless $r=2$, in which case $w$
visits the unique edge of the bracelet exactly $4$ times.

\subsubsection{} As in 
\S\ref{sec-almost2},
let $\check{w}$ be the result of dropping the last letter
of $w$ and let $\sigma$ be a cyclic permutation of
$\{1,\dots,\ell(\check{w})\}$. Then $\check{w}^\sigma \alpha_{\sigma(1)}$
is a critical weak Wigner word with graph,  bracelet and circuit length
the same as for $w$. We say that $\sigma$ is a {\em polarization} of $w$
if the last edge of $G$ visited by the walk
$\check{w}^\sigma \alpha_{\sigma(1)}$ is an edge of $Z$. 
Clearly:
\begin{itemize}
\item There exist exactly $2r$ polarizations of $w$. 
\end{itemize}
Note that the set of
polarizations of $w$ depends only on the equivalence class of $w$.

\subsubsection{}
Suppose now that we are given a polarization $\sigma$ of $w$.
We define the {\em canonical
decomposition}
$$\check{w}^\sigma=p_1p_2\cdots p_{2r-1}p_{2r}$$
associated to $\sigma$ to be the
unique decomposition with breaks at visits of the walk $\check{w}^\sigma$
to edges of the bracelet. From the bracelet-and-pendant-trees picture
it is
not difficult to deduce that each $p_i$ is a Wigner word and that no
two of the $p_i$ have letters in common with the exception that first
letters may coincide. Let
$s=\alpha_1\cdots
\alpha_{2r}$ be the sequence of first letters of the $p_i$. We call $s$
the {\em signature} associated to the critical weak Wigner word $w$ and
its polarization
$\sigma$. Necessarily
$s\alpha_1$ is a walk on the bracelet of $w$ visiting every edge of the
bracelet exactly twice unless $r=2$, in which case $s\alpha_1$ visits
the unique edge of the bracelet exactly $4$ times. Up to
equivalence of words there are very few possibilities for $s$. In fact,
the following possibilities are mutually exclusive and exhaustive:
\begin{itemize}
\item $r\geq 3$ and $s\sim 123\cdots r123\cdots r$.
\item $r=1$ and $s\sim 11$.
\item $r=2$ and $s\sim 1212$.
\item $r\geq 3$ and $s^\tau\sim 123\cdots r1r\cdots 2$
for some cyclic permutation $\tau$ of $\{1,\dots,2r\}$.
\end{itemize}
In the first case we say that the signature is {\em unidirectional},
whereas in the remaining cases we say that the signature is {\em
backtracking}.  Notice that if $s$ is unidirectional (resp.,
backtracking)
for some polarization $\sigma$, then $s$ is unidirectional (resp.,
backtracking) for all polarizations $\sigma$. Thus it makes sense to
say that $w$ itself is either unidirectional or backtracking.
\subsubsection{}
If $w$ is backtracking, then for some polarization
$\sigma$ the associated signature is of the form
$11$ if $r=1$, $1212$ if $r=2$, or $123\cdots r 1r\cdots 2$
if $r\geq 3$, in which case we say that $\sigma$ is a {\em strong
polarization} of $w$. It is not difficult to verify that:
\begin{itemize}
\item If $w$ is backtracking, there exist
exactly
$2$ strong polarizations of $w$ unless $r=2$, in which case
every polarization is strong (and so there exist exactly $4$
strong polarizations).
\end{itemize}
 Note that the set of strong polarizations of $w$ depends
only on the equivalence class of
$w$.

\subsection{Enumeration of critical weak Wigner words by Wigner words}
We again use the enumerative apparatus introduced in 
\S\ref{subsubsection:EnumerativeApparatus}.
We have
\begin{equation}\label{equation:CriticalWignerEnumeration}
\begin{array}{cl}
&\displaystyle\sum_{w} \varphi(w)\\\\
=&\displaystyle \sum_{r=1}^\infty
\sum_{u_1\in U_1}\cdots \sum_{u_{r}\in U_{r}}
\sum_{v_1\in V_1}\cdots \sum_{v_r\in V_r}\sum_\sigma\varphi(u^\sigma
\alpha_{\sigma(1)})\bigg/\left\{\begin{array}{cl} 4&\mbox{if  $r=2$}\\
2&\mbox{if $r\neq 2$}
\end{array}\right.\\\\
&\displaystyle+\sum_{r=3}^\infty
\sum_{x_1\in U_1}
\cdots \sum_{x_r\in U_r}
\sum_{y_1\in V_1}
\cdots \sum_{y_r\in V_r}
\sum_\tau \varphi(v^\tau \beta_{\tau(1)})/2r\\
\end{array}
\end{equation}
where:
\begin{itemize}
\item $w$ ranges over any cross-section of the set of critical weak
Wigner words;
\item $u=u_1\cdots u_{r}v_1v_r\cdots v_2=[\alpha_i]_{i=1}^{\ell(u)}$;
\item $\sigma$ ranges over cyclic permutations of $\{1,\dots,\ell(u)\}$;
\item $v=x_1\cdots x_ry_1\cdots y_r=[\beta_i]_{i=1}^{\ell(v)}$; and
\item $\tau$ ranges over cyclic permutations of $\{1,\dots,\ell(v)\}$.
\end{itemize}
In this setting we have
\begin{equation}\label{equation:CriticalWignerEnumerationM}
M(u^\sigma\alpha_{\sigma(1)})=
\prod_{i=1}^r(M(u_i)M(v_i))\cdot \left\{\begin{array}{ll}
d^{(2)}(\kappa(\gamma_1))&\mbox{if $r=1$,}\\
s^{(4)}(\kappa(\gamma_1),\kappa(\gamma_2))&\mbox{if $r=2$,}\\
K_r(\kappa(\gamma_1),\dots,\kappa(\gamma_r))&\mbox{if $r\geq 3$,}
\end{array}\right.
\end{equation}
and we have an analogous expression for $M(v^\tau\beta_{\tau(1)})$.
Formula (\ref{equation:CriticalWignerEnumeration}) may be derived from the
preceding discussion of
the bracelet of a critical weak Wigner word in a straightforward way.
 We omit the details.

\subsection{End of the proof}
The left sides of (\ref{equation:MeanFormula})
and (\ref{equation:CLTScaleBis}) coincide by
 formula (\ref{equation:WordExpansion5})
coming up in the proof of Theorem~\ref{Theorem:Preliminary}.
So we can rewrite (\ref{equation:CLTScaleBis})
as an identity of
formal power series
\begin{equation}\label{equation:CLTScaleTer}
\begin{array}{cl}
&\displaystyle
\sum_{n=1}^\infty 
\left(\lim_{k\rightarrow\infty}
 N_k\cdot\left(\langle \overline{L}(\NNN_k),x^n\rangle-
\langle \mu_{\NNN_k},x^n\rangle
\right)\right)t^{n+1}\\\\
=&\displaystyle
-\frac{1}{2}\sum_{[w,\alpha,\beta]}EM(w,[\beta\mapsto
\alpha])H([\beta\mapsto \alpha]_*w,t)+\sum_u
EMH(u,t)
\end{array}
\end{equation}
where:
\begin{itemize}
\item $[w,\alpha,\beta]$ ranges over a cross-section of the set of marked
Wigner words; and
\item $u$ ranges
over a cross-section of the set of critical weak
Wigner words.
\end{itemize}
To finish the proof of the theorem we have just to make the right side of
(\ref{equation:CLTScaleTer}) explicit. This can be done by exploiting
(\ref{equation:MarkedWignerEnumeration}),
(\ref{equation:MarkedWignerEnumerationM}),
(\ref{equation:CriticalWignerEnumeration}), and
(\ref{equation:CriticalWignerEnumerationM}).
Note that many of the terms in
the sum on $[w,\alpha,\beta]$ are cancelled by terms in the sum on
$u$ due to the parallel structure of formulas
(\ref{equation:MarkedWignerEnumerationM}) and
(\ref{equation:CriticalWignerEnumerationM}).
We
omit the remaining details of the proof because the calculations  are
very similar to those undertaken to prove
Lemma~\ref{equation:LemmaKeyCovariance}. The proof of
Theorem~\ref{Theorem:MainBis} is complete.
\qed

\section{Concentration}
\label{sec-poincareconc}
In this section
we work out sufficient conditions allowing one
to prove a CLT
for test functions more general than polynomials. Toward this end,
we  define for random matrices a notion of concentration 
and  a notion of CLT for polynomial test functions.
Then, assuming concentration, a polynomial-type CLT,
and a further
condition on the limiting covariance for polynomial test functions, we
prove a CLT for continuously differentiable test functions with
polynomial growth (Proposition~\ref{Theorem:generalCLT}).  Furthermore,
we establish the concentration property for the matrices $X(\NNN_k)$ 
studied in  Theorems 
\ref{Theorem:Preliminary} and \ref{Theorem:Main}
when the random variables $\xi_{\{\alpha,\beta\}}$
satisfy the \Poincare inequality 
with the same constant (Proposition~\ref{theo-poincaconc}). The main
result of this section (Theorem
\ref{cor-cltpoincare}) summarizes
the preceding considerations in a fashion convenient for applications in
\S\ref{sec-ChebDiag}.

\subsection{The concentration property}
\label{subsec-11.1}
Throughout this section
$\{Y_k\}_{k=1}^\infty$ 
denotes a sequence of random symmetric matrices.  
For such a general sequence we
are going to define and study a concentration property.
Eventually we are going to take
$Y_k=X(\NNN_k)$,
but in anticipation of applications of the concentration idea beyond the
scope of this paper, we work in a general setting until the end of the
proof of Proposition~\ref{Theorem:generalCLT}.
For any Lipshitz function $g:\RR^n\rightarrow\RR$
set
\begin{equation}\label{equation:LipDef}
\|g\|_{\rm Lip}:=
\sup_{
\begin{subarray}{c}
x, y\in \RR^n\\
x\neq y
\end{subarray}} \frac{|g(x)-g(y)|}{|x-y|}\,,
\end{equation}
where $|x-y|$ is the Euclidean distance between $x$ and $y$.
\begin{Definition}
We say that the sequence of matrices
$\{Y_k\}_{k=1}^\infty$
satisfies the {\em concentration
property} under the following conditions: 
\begin{equation}
\label{eq-150304a}
\begin{array}{l}
\mbox{\it There exists a 
constant $c>0$ such that
for any Lipschitz function}\\
\mbox{\it $g:\RR\rightarrow\RR$, it holds that}\
\sup_k\Var \trace g(Y_k)\leq c \|g\|_{\rm Lip}^2\,.
\end{array}\end{equation}
\begin{equation}
\label{eq-150304b}
\begin{array}{l}
\mbox{\it There exists a compact set $S\subset\RR$ such
that for any function}\\
\mbox{\it $f:\RR\rightarrow\RR$ supported in $S^c$ of polynomial
growth, it holds that}\\
E\left([\trace f(Y_k)]^2\right)\to_{k\to\infty} 0.
\end{array}
\end{equation}
\end{Definition}
\noindent

\noindent The next lemma deduces from
(\ref{eq-150304a}) and (\ref{eq-150304b})  a single statement
convenient for applications:
\begin{Lemma}\label{Lemma:Protection}
Suppose that $\{Y_k\}_{k=1}^\infty$ 
satisfies the concentration property. 
Then there exists a constant
$\bar{c}>0$
and a compact interval $T$
such that for any function
$f$ continuously differentiable on $T$
and of polynomial growth one has 
$$\limsup_{k\rightarrow\infty}
\Var \trace f(Y_k)\leq \bar{c}\sup_{x\in T}|f'(x)|^2.$$
\end{Lemma}
\proof Let $S$
be as in
(\ref{eq-150304b}). Choose a compact
interval
$I$ with interior containing the set
$S$,
and then choose a compact interval
$T$ with interior containing
$I$. Let
$g:\RR\rightarrow [0,1]$ be a continuously differentiable function
identically equal to
$1$ on
$I$ and identically vanishing in the complement of $T$. Let $\ell$ be
the length of $T$. Without loss of generality we may assume that
$f$ vanishes at some point of $T$. Then
$$\|fg\|_{\rm Lip}\leq \left(1+\ell\sup_{t\in T}|g'(t)|\right)\sup_{t\in
T}|f'(t)|,\;\;\;
\supp f(1-g)\subset S^c,$$
and
$$\left[\Var \trace f(Y_k)\right]^{1/2}
\leq \left[\Var\trace (fg)(Y_k)\right]^{1/2}+
(E[\trace (f(1-g))(Y_k)]^2)^{1/2},$$
whence the result by definition of the concentration property.
\qed 
\subsection{CLT's for differentiable test functions}
\label{subsec-11.2}
Our goal is to prove under suitable hypotheses
a central limit theorem for random variables of the form
$$Z_{f,k}:=\trace f(Y_k)-E\trace f(Y_k)$$
where $f$ is continuously differentiable on a large enough compact set
and of polynomial growth.
\begin{Definition}
We say that the sequence $\{Y_k\}_{k=1}^\infty$ satisfies a {\em
polynomial-type CLT} if there exists a mean zero Gaussian family
$\{W_n\}_{n=0}^\infty$ of random variables such that for every
polynomial function $f(x)=\sum_{i=0}^{m} a_ix^i$ it holds that
$Z_{f,k}$ converges in distribution as
$k\rightarrow\infty$ to $W_f:=\sum_{i=0}^{m} a_i W_i$.
\end{Definition}
\noindent The next proposition gives hypotheses under which one can
extend a CLT statement from polynomial test functions to differentiable
test functions of polynomial growth.
After proving the proposition, verification of its
 hypotheses for 
$Y_k=X(\NNN_k)$ under the assumptions of Theorems
\ref{Theorem:Preliminary} and
\ref{Theorem:Main}, along with further structural assumptions
concerning the functions $d^{(2)}$, $s^{(2)}$ and $s^{(4)}$
will be our task for the rest of the paper.

\begin{Proposition}\label{Theorem:generalCLT}  
Assume that the sequence of matrices $\{Y_k\}_{k=1}^\infty$ satisfies
both the concentration property and a polynomial-type CLT.
%
Assume 
further the existence 
of a sequence $\{q_n\}_{n=1}^\infty$ of polynomial functions with
the following properties:
\begin{itemize} 
\item For some compactly supported finite measure $\nu$ on $\RR$
the sequence $\{q_n\}_{n=1}^\infty$
is an orthonormal system in
$L^2(\nu)$. 
\item Every polynomial in $x$ is a finite linear combination
of the  $q_n(x)$. 
\item With $\bar q_n(x):=\int_0^x q_n(y) dy$,
the covariance matrix $K(m,n):=EW_{\bar q_m}W_{\bar
q_n}$ of the mean zero Gaussian family $\{W_{\bar q_n}\}_{n=1}^\infty
$
is diagonal.
\end{itemize}
Fix $T$ and $\bar{c}$ as in Lemma~\ref{Lemma:Protection}, with
$T\supset \supp\nu$.  Then, for any 
function $f$ of polynomial growth
which is continuously  differentiable on $T$,
the random variables $Z_{f,k}$
converge in distribution to  a  mean  zero
Gaussian random variable $Z_f$ with variance
\begin{equation}
\label{eq:300304f}
EZ_f^2=\|f'\|_K^2\leq \bar{c}\,\sup_{t\in T}|f'(t)|^2,
\end{equation}
where for any function $h$ continuous on $T$ we set
$$\|h\|_K^2:=\sum_{n=1}^\infty K(n,n) \langle \nu,hq_n\rangle^2.$$
\end{Proposition}
\proof Consider at first the case in which $f$ is a polynomial.
The polynomial-type CLT implies that the variables
$Z_{f,k}$ converge  in distribution to
$W_f$,
and since $f$ differs by a constant from a finite linear combination of
the
$\bar{q}_n$, the variance $EW_f^2$ takes
the value asserted in 
(\ref{eq:300304f}), namely $\|f'\|^2_K$. Furthermore, 
by 
Lemma~\ref{Lemma:Protection} and the Fatou Lemma, the estimate for
$\|f'\|_K^2$ asserted in (\ref{eq:300304f}) holds. Thus all assertions
are proved if
$f$ is a polynomial function.

We turn to consideration of the general case.
Let $\{Q_m\}_{m=1}^\infty$ be a sequence of polynomials 
tending uniformly on $T$ to $f'$ (such is provided by the
Stone-Weierstrass theorem) and put
$\overline{Q}_m(x):=\int_0^xQ_m(y)dy$.
Clearly, the sequence
$\{Q_m\}_{m=1}^\infty$ is
$\|\cdot\|_{L^2(\nu)}$-Cauchy. But by (\ref{eq:300304f}) the sequence
$\{Q_m\}_{m=1}^\infty$
is also
$\|\cdot\|_K$-Cauchy. A dominated convergence
argument now shows that
$\|f'\|_K^2=\lim_{m\rightarrow\infty}\|Q_m\|_K^2$.
It follows
that the estimate for
$\|f'\|_K^2$ asserted in
(\ref{eq:300304f}) holds. 
By Lemma~\ref{Lemma:Protection} the family of random variables $Z_{f,k}$
is tight; let $Y$ be any subsequential limit-in-distribution. 
For any $t\in \RR$ one has
$$|E e^{itY}-Ee^{itW_{\overline{Q}_m}}|\leq
\limsup_{k\rightarrow\infty}
E|e^{itZ_{f-\overline{Q}_m,k}}-1|
\leq |t|\limsup_{k\rightarrow\infty}(EZ_{f-\overline{Q}_m,k}^2)^{1/2}.$$
The
quantity on the right by
Lemma~\ref{Lemma:Protection} tends to $0$ as
$m\rightarrow\infty$, and clearly
$$Ee^{itW_{\overline{Q}_m}}=
e^{-t^2\|Q_m\|_K^2/2}\rightarrow_{m\rightarrow\infty}
e^{-t^2\|f'\|_K^2/2}.$$ Therefore (the characteristic function of)
$Y$ is (that of)
a mean zero Gaussian random variable of variance $\|f'\|_K^2$. Since all
subsequential limits are the same we get convergence-in-distribution of
$Z_{f,k}$ to a mean zero Gaussian random variable of variance
$\|f'\|^2_K$.  All
assertions have been proved.
\qed 

\subsection{\Poincare inequalities for matrices}
\label{subsec-poincarematrices}
We say that a probability
distribution $\eta$ on $\RR$
satisfies a {\em \Poincare
inequality} if there exists a constant $c_\eta$ such that for
any $f:\RR\rightarrow \RR$ smooth, it holds that
$$\Var_\eta(f):= \int \left(f(x)-\int f(x)\eta(dx)\right)^2 \eta(dx)\leq 
c_\eta \int |f'(x)|^2 \eta(dx)\,.$$
For such a distribution $\eta$ one has
\begin{equation}
\label{eq-290304a}
E \exp \left(\frac{|Y-E Y|}{12\sqrt{c_\eta}}\right)\leq
2\,\;\;\;(Y:\mbox{ random variable with law $\eta$}),
\end{equation}
see \cite[Theorem 2]{borovkov-utev} (or \cite{bobkov} for optimal
constants).

It is well known (see, e.g., \cite[Pg. 49]{Ledo}) that if $\eta_i,
i=1,\ldots,K$ satisfy
\Poincare inequalities with constants $c_{\eta_i}$, then for any
smooth function $g:\RR^K\rightarrow \RR$, and with $\eta=
\otimes_{i=1}^K \eta_i$ and
$c_\eta=\max_{i=1}^K c_{\eta_i}$, one has 
\begin{equation}
\label{eq-poincare}
\Var_\eta(g)=: \int \left(g(x)-\int g(x)\eta(dx)\right)^2 \eta(dx)\leq 
c_\eta \int |\nabla g(x)|^2 \eta(dx)\,.
\end{equation}
We recall (see e.g. \cite[Lemma 1.2]{GZconc}) that if $f:\RR\rightarrow
\RR$  is Lipschitz with Lipschitz constant $\|f\|_{\rm Lip}$, 
then the function $f_N:\RR^{ N( N+1)/2}\rightarrow \RR$ 
on $N$-by-$N$ symmetric matrices given by 
$f_N(X)=\trace f(X)$ 
is Lipschitz with  Lipschitz constant 
$\|f_N\|_{\rm Lip}\leq \sqrt{ N} \|f\|_{\rm Lip}$.
It follows that if $X$ is an $N$-by-$N$ symmetric random matrix
with on-or-above-diagonal entries independent and
satisfying the \Poincare inequality with the same constant $c/ N$,
then  for any Lipshitz $f:\RR\rightarrow\RR$ one has
$$\Var \trace f(X)\leq c \|f\|_{\rm Lip}^2\,.$$
See
\cite{Chabo} for a systematic use of this fact, and 
\cite{GZconc} for other concentration inequalities for random 
matrices. In particular, in the setting of Theorems
\ref{Theorem:Preliminary} and \ref{Theorem:Main},
if the random variables $\xi_{\{\alpha,\beta\}}$
 satisfy the Poincar\'{e} inequality with the same
constant $c$, we see that  (\ref{eq-150304a}) above holds true,
with $Y_k=X(\NNN_k)$.

\begin{Proposition}
\label{theo-poincaconc}
In the setting and under the hypotheses of Theorems
\ref{Theorem:Preliminary} and \ref{Theorem:Main},
suppose that the random variables $\xi_{\{\alpha,\beta\}}$
satisfy the \Poincare inequality with the same constant $c$. Then
the sequence $\{X(\NNN_k)\}_{k=1}^\infty$
has the concentration property.
\end{Proposition}

\noindent Before proving the proposition, we state an auxiliary
estimate, which may be of interest in its own right. We get
the estimate by combining the ideas of \cite{furedi} as
summarized in Proposition~\ref{Proposition:FK} above with the moment
bound (\ref{eq-290304a}).  We remark in passing that under somewhat
stronger assumptions, a considerably stronger assertion could be obtained
by the methods of
\cite{soshnikov}.

\begin{Lemma}
\label{lem-expmompoi}
Under the assumptions of Proposition \ref{theo-poincaconc}, there
exist constants \linebreak $C>0$ and $\epsilon>0$ such that with
$r(N):=\lfloor N^\epsilon\rfloor$ one has
\begin{equation}
\label{eq-190304a}
\frac1N E\trace X(\NNN_k)^{2r(N_k)}\leq C^{2r(N_k)} \,
\end{equation}
for all sufficiently large $k$.
\end{Lemma}
\proof By an obvious rescaling, we may assume without loss of
generality that
$C(2)=1$, where $C(2)$ is as defined in (\ref{eq:300304e}).
We may further assume, so we claim, that $D=0$. 
To see that this is so, set $(\overline X(\NNN_k))_{\alpha\beta}=
N_k^{-1/2}\xi_{\{\alpha,\beta\}}$ and suppose that the lemma holds
with $\overline X(\NNN_k)$ in place of $X(\NNN_k)$.
Then 
\begin{eqnarray*}
E\trace X(\NNN_k)^{2r(N_k)}&\leq&
N_k |\lambda|_{\max}(X(\NNN_k))^{2r(N_k)}\leq 
 N_k[|\lambda|_{\max}(\overline X(\NNN_k))+|D|_\infty]^{2r(N_k)}\\
&\leq&
2^{2r(N_k)-1} 
N_k [|\lambda|_{\max}(\overline X(\NNN_k))^{2r(N_k)}+
|D|_\infty^{2r(N_k)}]
\\&\leq& 2^{2r(N_k)-1}N_k
[E\trace \overline X(\NNN_k)^{2r(N_k)}+|D|_\infty^{2r(N_k)}]\\
&
\leq & 2^{2r(N_k)-1}N_k^2
[C^{2r(N_k)}+ |D|_\infty^{2r(N_k)}]\leq (3(C+|D|_\infty))^{2r(N_k)}\,,
\end{eqnarray*}
for  all $k$ large enough. The claim is proved.
We assume for the rest of the proof that $D=0$. 
By (\ref{equation:WordExpansion2}), for any positive integer $n$,
\begin{equation}
\label{eq:290304c}
\langle \overline{L}(\NNN_k),x^{2n}\rangle\leq
\sum_{q=1}^{n+1} N_{\rm FK}(2n,q) 
 N_k^{q-(n+1)} \max_{b\in {\rm FK}(2n,q)}
E|\xi(b)|
\end{equation}
with ${\rm FK}(2n,q)$ denoting the collection of 
weak Wigner words of length $2n+1$ and
weight $q$, and $N_{\rm FK}$ as in (\ref{eq:290304b}). 
Note that
\begin{eqnarray*}
&& \max_{b\in {\rm FK}(2n,q)}
E|\xi(b)| \leq  C(3(2n+2-2q))\\
&\leq& 
\sup_{\alpha,\beta}
E\left(\exp\left(|\xi(\alpha,\beta)|/12\sqrt{c}\right)\right)
(1\vee (12\sqrt{c}))^{3(2n+2-2q)}[3(2n+2-2q)]!\\
&\leq&
2(1\vee (12\sqrt{c}))^{3(2n+2-2q)}[3(2n+2-2q)]!=:
2C_0^{3(2n+2-2q)}[3(2n+2-2q)]!\,,
\end{eqnarray*}
where the first inequality is due to
(\ref{eq:300304e})
and the second to (\ref{eq-290304a}).
Thus
\begin{eqnarray*}
\langle \overline{L}(\NNN_k),x^{2n}\rangle&\leq&
2^{n+1}\sum_{q=1}^{n+1} 
N_k^{q-(n+1)} {[3(2n+2-2q)]! (C_0n)^{3(2n+2-2q)}}\\
&\leq&
2^{n+1}\sum_{j=0}^n N_k^{-j} (6C_0n)^{12j}
\leq 2^{n+2}
\end{eqnarray*}
as long as $(6C_0n)^{12}/ N_k\leq 1/2$. This completes the proof.
\qed\\

\noindent
{\it Proof of Proposition \ref{theo-poincaconc}.}\,
In view of Theorem \ref{Theorem:Preliminary} 
and the discussion in \S\ref{subsec-poincarematrices},
it only
remains to check (\ref{eq-150304b}). This is based on
Lemma
\ref{lem-expmompoi}. Fix $C$ as in the statement of that
lemma.
Define the compact set $S=[-C-1,C+1]$. Suppose that $|f(x)|\leq c_1
|x|^{c_2}$ and $f$ is supported on $S^c$. Then, using that
$$(|x|/(C+1/2))^{r(N_k)}\geq |x|^{2c_2}\;\;\;\mbox{for $|x|\geq C+1$
and $k$ large},$$
one has
\begin{eqnarray}
\label{eq-190304b}
E\left([\trace f(X(\NNN_k))]^2\right)&\leq &
 N_k  
E\trace f^2(X(\NNN_k))\nonumber \\
&\leq &
 N_k c_1^2  
E\sum_{i=1}^{ N_k} 
\lambda_i(\NNN_k)^{2c_2} {\bf 1}_{|\lambda_i(\NNN_k)|
\geq (C+1)}
\nonumber \\
&\leq &
 N_k c_1^2  
E\sum_{i=1}^{ N_k} \left(\frac{\lambda_i(\NNN_k)}{C+1/2}
\right)^{r(N_k)}
\nonumber \\
&\leq &
   N_k^2c_1^2 
\left(\frac{C}{C+1/2}\right)^{r(N_k)}\to_{k\to\infty } 0\,.
\end{eqnarray}
\qed

By combining Propositions \ref{Theorem:generalCLT}
and \ref{theo-poincaconc},
we immediately get the following theorem, which is
the main result of this section.
Recall that under the assumptions of Theorem \ref{Theorem:Main}, 
the sequence $\{X(\NNN_k)\}_{k=1}^\infty$ satisfies a polynomial-type
CLT, i.~e., there exists a mean zero Gaussian family
$\{W_n\}_{n=0}^\infty$ of random variables such that for every polynomial
function $f(x)=\sum_{i=0}^m a_ix^i$ the random variables $\trace
f(X(\NNN_k))-E\trace f(X(\NNN_k))$ converge in distribution as
$k\rightarrow\infty$ to $W_f:=\sum_{i=0}^m a_iW_i$.

\begin{Theorem}
\label{cor-cltpoincare}
We work in the setting and under the hypotheses of Theorems
\ref{Theorem:Preliminary} and \ref{Theorem:Main}. 
We make the following further assumptions:
\begin{itemize}
\item The random variables $\xi_{\{\alpha,\beta\}}$ 
satisfy the Poincar\'{e} inequality with the same constant $c$
(and hence $\{X(\NNN_k)\}_{k=1}^\infty$
has the concentration property). 
\item There exists a sequence $\{q_n\}_{n=1}^\infty$ of polynomial
functions with the following properties:
\begin{itemize}
\item For some compactly supported finite measure $\nu$ on
$\RR$, the sequence $\{q_n\}_{n=1}^\infty$ is an orthonormal system in
$L^2(\nu)$.
\item Every polynomial in $x$ is a finite linear combination
of the  $q_n(x)$.
\item With $\bar q_n(x):=\int_0^x q_n(y) dy$, the
covariance matrix
$K(m,n):=EW_{\bar q_m}W_{\bar q_n}$ of the mean zero Gaussian family
$\{W_{\bar q_n}\}_{n=1}^\infty$ is diagonal.
%
\end{itemize}
\end{itemize}
Then there exists a compact interval $T\supset\supp\nu$ and a constant 
$\bar{c}>0$
such that for 
any
function $f$ of polynomial growth
which is continuously  differentiable on $T$,
the random variables 
$$Z_{f,k}:= \trace f(X(\NNN_k))-E\trace f(X(\NNN_k))$$
converge in distribution to  a  mean  zero
Gaussian random variable $Z_f$ with variance
\begin{equation}
\label{eq:300304fnew}
EZ_f^2=
\sum_{n=1}^\infty K(n,n) \langle \nu,f'q_n\rangle^2
\leq \bar{c}\,\sup_{t\in T}|f'(t)|^2\,.
\end{equation}
\end{Theorem}

%
            
%

\section{Diagonalization by Chebyshev polynomials}
\label{sec-ChebDiag}
We discuss two
specializations of the band matrix model in which we can make
$\mu$, $\Phi(c,t)$, $\Theta(x,y)$, $\Psi(x,y)$,
 $\Var Z_f$ and $E_f$ as appearing in Theorems~\ref{Theorem:Preliminary}
and \ref{Theorem:Main} much
more explicit and moreover apply Theorem~\ref{cor-cltpoincare}.
This will be possible because in these specializations (slight variants
of) Chebyshev polynomials  diagonalize the covariance matrix of the
limiting mean zero Gaussian random variables.

\subsection{Inversion of power series and $p$-Chebyshev
polynomials} Our computation involves the inversion of  formal power
series.  Fix a sequence of real numbers $\{a_i\}$ and define the formal
power  series
\begin{equation}
\label{eq:040604b}
p=p(t):=t+\sum_{i=2}^\infty a_it^i\,.
\end{equation}
(For the proof of Theorem~\ref{Theorem:Mainspecial} concerning the
generalized Wigner matrix model, it will be enough simply to take
$p(t)=\Phi(t)$, where
$\Phi(t)$ is the generating function for the Catalan numbers defined in
(\ref{eq:040604a}) below.) For each
positive integer
$n$, define the {\em $n^{th}$ $p$-Chebyshev polynomial}  
$T_{n,p}(x)$
as the unique polynomial in $x$ of degree $n$
with real coefficients such that $T_{n,p}(1/t)$ is the principal part of
the Laurent
series $p(t)^{-n}$. Finally,
define the matrix $P$ with rows and columns indexed by
the positive integers  by setting $P_{ij}$ equal to the coefficient of
$t^j$ in $p^i$, i.~e.,
\begin{equation}
\label{eq:040604c}
P_{ij}:=\Res_{t=0}\left(t^{-j}p^i\frac{dt}{t}\right)\,,
\end{equation}
where for any sequence $[c_i]_{i=-\infty}^\infty$ of constants
such that $c_i=0$ for
$i\ll 0$ we set
$$\Res_{t=0}\sum_{i=-\infty}^\infty c_it^i\,dt:=c_{-1}\,.$$
\begin{Lemma}\label{Lemma:Lagrange}
Fix $p(t)$ as in (\ref{eq:040604b}) with its associated $p$-Chebyshev 
polynomials $T_{n,p}(x)$
and matrix $P$ as in (\ref{eq:040604c}).
Identify power series in $x$ without constant term
in the obvious way with column vectors having entries indexed
by the positive integers (thus identifying polynomials in $x$ without
constant term with finitely supported infinite column vectors).
Then, the
$n^{th}$ column of
$P^{-1}$ equals $\frac{1}{n}xT_{n,p}'(x)$.
\end{Lemma}

\proof Let $r=r(t)$ be the formal power series inverse of $p(t)$,
i.~e., the unique power series without constant term such that
$$p(r(t))=r(p(t))=t.$$
By the  {\em Lagrange inversion
formula}, c.f. \cite[\S 5.4]{stanley}, 
$$ (P^{-1})_{ij}
=\Res_{t=0}\left(t^{-j}r^i\frac{dt}{t}\right)=
\frac{i}{j}
 \Res_{t=0}\left(p^{-j}t^i\frac{dt}{t}\right).$$
The last expression is by definition exactly the coefficient of
$x^i$ in $\frac{1}{j} x T_{j,p}'(x)$.
\qed

\subsection{Chebyshev polynomials}

\subsubsection{Definition}
For each positive integer $n$ we define the $n^{th}$
{\em Chebyshev polynomial} $T_n(x)$ of the {\em first kind} to be
the unique polynomial in $x$ such that
$$T_n(z+1/z)=z^n+1/z^n\;\;\;
(\mbox{equivalently:}\;T_n(2\cos \theta)=2\cos n\theta)$$
and we define the
$n^{th}$ {\em Chebyshev polynomial}
$U_n(x)$ of the {\em second kind} by the rule
$$U_n(x):=\frac{1}{n}T_n'(x).$$
We have orthogonality relations
\begin{equation}\label{equation:SemicircleOrthogonality}
\frac{1}{2\pi}\int_{-2}^2U_m(x)U_n(x)\sqrt{4-x^2}\,dx
=\delta_{mn}\;\;\;(m,n=1,2,3,\dots)
\end{equation}
as can be
verified directly by the trigonometric substitution
$x=2\cos \theta$.  Note that the weight
figuring in these orthogonality relations is the semicircle law
$\sigma_S$ of mean
$0$ and variance $1$. These relations say that the family
$\{U_n(x)\}_{n=1}^\infty$ is the Gram-Schmidt orthogonalization in
$L^2(\sigma_S)$ of the family $\{x^{n-1}\}_{n=1}^\infty$ of powers of
$x$. We have analogous orthogonality relations for Chebyshev polynomials
of the first kind (with a different weight), but these we omit because we
have no use for them.
\subsubsection{Warning}
Our definitions are
not quite the standard ones. One usually defines
$T_n(\cos
\theta)=\cos n\theta$ and $U_n(x)=\frac{1}{n+1}T_{n+1}'(x)$. We have
rescaled and re-indexed in order to obviate many annoying factors of $2$
and shifts of $1$. It is also worth pointing out that in our set up 
the polynomial
$xU_n(x)$ is monic of degree $n$, and moreover even or odd
according as $n$ is even or odd.

\subsubsection{Reinterpretation of the Chebyshev polynomials}
\label{subsubsection:PhiChebyshev}
Consider the odd power series
\begin{equation}
\label{eq:040604a}
\Phi(t)=\frac{1-\sqrt{1-4t^2}}{2t}=
\sum_{n=0}^\infty\frac{1}{n+1}\left(\begin{array}{c}2n\\
n\end{array}\right)t^{2n+1}= t+t^3+2t^5+5t^7+\cdots
\end{equation}
having the $n^{th}$ Catalan number as the coefficient of
$t^{2n+1}$.
Clearly $\Phi(t)$ satisfies 
the functional equation
$$
1/t=\Phi(t)+1/\Phi(t)
$$
 and hence more generally the functional equation
\begin{equation}
\label{eq:040604e}
T_n(1/t)=\Phi(t)^n+1/\Phi(t)^n
\end{equation}
for all positive integers $n$. 
In other words, for each positive integer $n$,
the $n^{th}$ Chebyshev polynomial $T_n(x)$ (with its constant term dropped)
may be reinterpreted as the
$n^{th}$
$\Phi$-Chebyshev polynomial $T_{n,\Phi}(x)$, in the sense of
Lemma~\ref{Lemma:Lagrange}.

\subsubsection{Diagonalization formulas}
In
Lemma~\ref{Lemma:Lagrange} let us now
take $p(t)=\Phi\left(\frac{t}{1-\gamma t}\right)$
where $\gamma$ is any real constant. (For the proof of
Theorem~\ref{Theorem:Mainspecial} concerning the generalized Wigner
matrix model it will be enough to consider just the case
$\gamma=0$.) Note that 
$T_{n,p}(x)$ and $T_n(x-\gamma)$ differ by a constant and hence 
$\frac{1}{n}xT_{n,p}'(x)=xU_n(x-\gamma)$. Using the obvious
identification between power series in $x$ without constant term and
{\em row} vectors with entries indexed by the positive integers, one may
think of
$p^i(x)$ as $e_i^T P$, where $P$ is the matrix from
(\ref{eq:040604c}) and $e_i$ is the (infinite)  column vector whose
$j^{th}$ entry is
$\delta_{ij}$.   On the other hand, by the remark following
(\ref{eq:040604e}), and 
Lemma~\ref{Lemma:Lagrange}, using the obvious identification between
power series in $x$ without constant term and {\em column} vectors with
entries indexed by the positive integers, one can identify 
$xU_n(x-\gamma)$ with $P^{-1} e_n$. Hence,
for any sequence $\{\eta_i\}_{i=1}^\infty$ of real
constants, 
\begin{equation}\label{equation:ChebDiag1}
\left\langle \sum_{i=1}^\infty
\eta_i\Phi\left(\frac{t}{1-\gamma
t}\right)^{i},tU_{n}(t-\gamma)\right\rangle= \eta_n
\end{equation}
and similarly
\begin{equation}\label{equation:ChebDiag2}
\left\langle \sum_{i=1}^\infty \eta_i \Phi\left(\frac{x}{1-\gamma
x}\right)^i\Phi\left(\frac{y}{1-\gamma y}\right)^i,
xU_m(x-\gamma)yU_n(y-\gamma)\right\rangle=
\eta_n\delta_{mn},
\end{equation}
for all positive integers $m$ and $n$. 

\subsection{First specialization: generalized Wigner matrices}
\label{sec-firstgen}
\subsubsection{Specialization of the model}
As in Theorem~\ref{Theorem:Mainspecial}, assume now that 
$$D\equiv 0,\;\;\;\int s^{(2)}(c,c')\theta(dc')\equiv 1.
$$
This
specialization of the band matrix model we call the {\em generalized
Wigner matrix} model. In the case that $s^{(2)}\equiv 1$
this is more or less the standard Wigner matrix model, whence the
terminology.

\subsubsection{Calculation of $\Phi(c,t)$ and $\mu$}
In the case at hand $\Phi(c,t)$ must be independent of $c$,
and hence by functional equation (\ref{equation:PhiRecursion})
we have $\Phi(c,t)=\Phi(t)$ 
where the latter is as defined in (\ref{eq:040604a}). From $\Phi(t)$ we
can read the moments of $\mu$. We conclude that $\mu$ is the semicircle
law $\sigma_S$ of mean $0$ and variance $1$. 

\subsubsection{Calculation of $\Theta(x,y)$ and $\Psi(x,y)$}

Because $\Phi(c,t)=\Phi(t)$, the integrals
figuring in the definitions of $\Theta(x,y)$ and $\Psi(x,y)$ 
greatly simplify. We find that
\begin{equation}\label{equation:DiagonalCoefficients}
\Theta(x,y)=\sum_{r=1}^\infty\lambda_r
\Phi(x)^r\Phi(y)^r,\;\;\;
\Psi(x,y)=\sum_{r=1}^\infty\epsilon_r
\Phi(x)^r\Phi(y)^r
\end{equation}
where
$$\begin{array}{rcl}
\lambda_r&=&\displaystyle\frac{1}{r}\int\cdots \int
K_r(c_1,\cdots,c_r)
\theta(dc_1)\cdots \theta(dc_r),\\\\
\epsilon_r&=&\left\{\begin{array}{cl}
\displaystyle\int
(d^{(2)}(c)-2s^{(2)}(c,c))\theta(dc)&\mbox{if $r=1$,}\\
\displaystyle\frac{1}{2}\int\int
(s^{(4)}(c_1,c_2)-3s^{(2)}(c_1,c_2)^2)\theta(dc_1)\theta(dc_2)&\mbox{if
$r=2$,}\\ 0&\mbox{if $r\geq 3$.}
\end{array}\right.
\end{array}
$$
\subsubsection{Calculation of $\Var Z_f$ and $E_f$}
Using the orthogonality relations
(\ref{equation:SemicircleOrthogonality}), for any polynomial function
$f$ we can write
$$f'(x)=\sum_{n=1}^\infty U_n(x) \,Ef'(S)U_n(S)$$
where $S$ is a random variable with standard semicircular law
$\sigma_S$. (Only finitely many nonzero terms appear in the sum.)
Then, using the diagonalization formulas (\ref{equation:ChebDiag2}) (with 
$\eta_i=\lambda_i$) and
(\ref{equation:ChebDiag1}) (with $\eta_{2i}=\lambda_i$ and
$\eta_{2i+1}=0$),
we can in the present specialization of the band matrix model
rewrite (\ref{equation:CovarianceFormula}) 
and (\ref{equation:MeanFormula}) in the form
\begin{equation}\label{equation:CovarianceFormulaWignerRewrite}
\Var Z_f=\sum_{r=1}^\infty
(2\lambda_r+\epsilon_r)(Ef'(S)U_r(S))^2,
\end{equation} 
\begin{equation}
\label{eq:300304b}
E_f=\sum_{r=1}^\infty 
\frac{1}{2}\left(\lambda_r+\epsilon_r\right) Ef'(S)U_{2r}(S).
\end{equation}

\subsection{Proof of Theorem \ref{Theorem:Mainspecial}}
\label{subsec-proofofMainspecial}
In view of Theorem \ref{Theorem:Main}, the discussion in
\S\ref{sec-firstgen} immediately above,
and Theorem
\ref{cor-cltpoincare}, all that remains to be done is to
take in the latter Theorem $\nu=\sigma_S$ and $q_n(x)=U_n(x)$.
\qed

\subsection{Second specialization: generalized  Wishart
matrices}
\label{subsec-secondspec}
\subsubsection{Square-root generalized Wishart matrices}
\label{subsec-sqrtwish}
Assume that color space is decomposed as a disjoint union 
$A\cup B$ and that
$$0<\theta(A)\leq \theta(B).$$
Put
$$\alpha=\sqrt{\frac{\theta(B)}{\theta(A)}}\geq 1,\;\;\;
\beta=\sqrt{\frac{\theta(A)}{\theta(B)}}\leq
1,\;\;\;\gamma=\alpha+\beta\geq 2.$$ 
Assume that $s^{(2)}$ and $s^{(4)}$ vanish
identically
on $(A\times A)\cup(B\times B)$ and that
$$\int_A s^{(2)}(\cdot,c')\theta(dc')=\beta\one_B,\;\;\;
\int_B s^{(2)}(\cdot,c')\theta(dc')=\alpha\one_A.$$
Assume that $D\equiv 0$. Assume that $d^{(k)}\equiv 0$ for all $k>0$.
 This
specialization of the band matrix model we call the {\em generalized
square-root
Wishart matrix} model.

\subsubsection{Generalized Wishart matrices}
\label{subsec-genwishmat}
As we are about to see, the machinery we developed
is well-suited to deal with square-root generalized
Wishart matrices. In applications, however, one is often 
interested in a slight variant. Write 
$$\NNN=\NNN^A\cup\NNN^B$$
where $\NNN^A$ (resp., $\NNN^B$) is the set of letters in $\NNN$ with
color in
$A$ (resp.,
$B$). As usual let $N$, $N^A$ and $N^B$ denote the
corresponding cardinalities. By re-arranging coordinates,
$X(\NNN)$ can be written in the form
$$X(\NNN)=\left[\begin{array}{cc}
0& Y(\NNN)\\
Y^T(\NNN)& 0
\end{array}\right]\,,$$
where the matrix $Y(\NNN)$ has rows indexed by $\NNN^A$
and columns indexed by $\NNN^B$. 
In this situation, 
we call the symmetric random matrices
$$W(\NNN)=Y(\NNN)Y(\NNN)^T$$ 
(rows and columns indexed by $\NNN^A$) {\em generalized Wishart
matrices}, and we are interested in the empirical distribution of
their eigenvalues  
$\{\lambda_i(W(\NNN))\}_{i=1}^{N^A}$ (all non-negative):
$$L_W(\NNN):=\frac{1}{N^A}
\sum_{i=1}^{N^A} \delta_{\lambda_i(W(\NNN))}
\,.$$
In the case that $s^{(2)}$ 
is constant on $(A\times B)\cup (B\times
A)$, the spectrum of $W(\NNN)$ 
is the same as that of standard Wishart matrices, 
whence the terminology.
Note that 
for any function $g(t)$ on $\RR_+$ with $g(0)=0$, and
setting $\tilde{g}(t)=g(t^2)$, one has
\begin{equation}
\label{eq:300304g}
N\langle L(\NNN),\tilde{g}\rangle=
\trace \tilde{g}(X(\NNN))=2\trace
g(W(\NNN))=2N^A\langle L_W(\NNN),g\rangle.
\end{equation}
Hence, once results (either LLN or CLT) are derived for
$X(\NNN)$, it is a simple exercise in book-keeping to transform them to
statements about $W(\NNN)$.
 
\subsubsection{Calculation of $\Phi(c,t)$}
Under our additional assumptions in \S\ref{subsec-sqrtwish},
we can write
$$\Phi(\cdot, t)=\one_A\Phi_A(t)+\one_B\Phi_B(t)$$
where $\Phi_A(t)$ and $\Phi_B(t)$ are color-independent.
Functional equation (\ref{equation:PhiRecursion}) 
in the case at hand specializes to the functional
equation
$$
\one_A\Phi_A(t)+\one_B\Phi_B(t)=
t\left(1-\beta\one_B\Phi_A(t)-\alpha\one_A\Phi_B(t)\right)^{-1},
$$
which in turn can be rewritten as the pair of functional equations
$$\Phi_A(t)=t(1-t\alpha\Phi_B(t))^{-1},\;\;\;
\Phi_B(t)=t(1-t\beta\Phi_A(t))^{-1}.
$$
After a straightforward calculation with formal power series we find
that 
\begin{equation}\label{equation:WishartSTransform}
\begin{array}{rcl}
\displaystyle \int
\Phi(c,t)\theta(dc)&=&\theta(A)\Phi_A(t)+\theta(B)\Phi_B(t)\\
&=&\displaystyle
\frac{1-\sqrt{(1-\gamma t^2)^2-4t^4}}{\gamma t}
=t\left(1+\frac{2}{\gamma}\Phi\left(\frac{t^2}{1-\gamma
t^2}\right)\right),
\end{array}
\end{equation}
\begin{equation}\label{equation:WishartCLTGen}
\Phi_A(t)\Phi_B(t)=
\frac{1-\gamma t^2-
\sqrt{(1-\gamma t^2)^2-4t^4}}
{2t^2}=\Phi\left(\frac{t^2}{1-\gamma t^2}\right),
\end{equation}
where $\Phi(t)$ is as in \S\ref{subsubsection:PhiChebyshev}.

\subsubsection{Calculation of $\mu$} From
(\ref{equation:WishartSTransform}) we know the moments of the measure
$\mu$  and moreover we can compare these moments to those of the
semicircle distribution. We find that
\begin{equation}\label{equation:WishartMuFormula}
\langle \mu,f\rangle=\sqrt{1-\frac{4}{\gamma^2}}\,
f(0)+\frac{1}{\gamma\pi}\int_{|x^2-\gamma|<2}
\frac{f(x)\sqrt{4-(x^2-\gamma)^2}}{|x|}dx.
\end{equation}
Now $\mu$ is the weak limit in probability of the empirical distributions
$L(\NNN_k)$. To calculate the corresponding limit $\mu_W$ of the
empirical distributions
$L_W(\NNN_k)$ we use the ``bookkeeping principle'' (\ref{eq:300304g}) to
find that
\begin{equation}\label{equation:WishartMuFormulaBis}
\langle
\mu_W,f\rangle=\frac{\gamma+\sqrt{\gamma^2-4}}
{4\pi}\int_{\gamma-2}^{\gamma+2}
\frac{f(x)\sqrt{4-(x-\gamma)^2}}{x}dx.
\end{equation}
See e.g. 
\cite{pastur-marchenko} for the latter result in the case 
of Wishart matrices.
Note that if $\theta(A)<\theta(B)$ and hence $\gamma> 2$, the
measure
$\mu$ has some mass concentrated at the origin. Notice also that
if $\gamma=2$, then $\mu$ is the semicircle distribution.

\subsubsection{Calculation of $\Theta(x,y)$ and $\Psi(x,y)$}
With $\lambda_{2r}$ and $\epsilon_{2r}$ as 
defined in
(\ref{equation:DiagonalCoefficients}), 
we have
$$\begin{array}{rcl}
\Theta(x,y)&=&\displaystyle\sum_{r=1}^\infty
\lambda_{2r}\Phi\left(\frac{x^2}{1-\gamma x^2}\right)^r
\Phi\left(\frac{y^2}{1-\gamma y^2}\right)^r,\\\\
\Psi(x,y)&=&\displaystyle\sum_{r=1}^\infty
\epsilon_{2r}\Phi\left(\frac{x^2}{1-\gamma x^2}\right)^r
\Phi\left(\frac{y^2}{1-\gamma y^2}\right)^r.
\end{array}
$$
To verify these
formulas the main thing to note is that
$K_{r}(c_1,\dots,c_{r})=0$ unless
colors along the sequence $c_1,\ldots,c_r,c_1$
alternate between $A$ and $B$.
It follows in particular that
$K_r\equiv 0$ for odd $r$.

\subsubsection{The measure $\nu$ and associated orthogonal polynomials}
\label{subsubsection:NuDef}
Let $\nu$
be the measure with density 
$$\frac{d\nu}{dx}=\one_{|x^2-\gamma|\leq
2}\frac{\sqrt{4-(x^2-\gamma)^2}}{2\pi|x|}$$ with respect to Lebesgue
measure. Note that $\mu$ is a convex combination of $\nu$
and a unit mass at the origin. Note that if
$\gamma=2$ then
$\nu=\sigma_S$. 
Put $$V_n(x):=xU_n(x^2-\gamma).$$
By a straightforward calculation one verifies that the system of
polynomial functions $\{V_n\}_{n=1}^\infty$ is orthonormal in
$L^2(\nu)$, and moreover forms the ``odd part'' of the family of
orthogonal polynomials naturally associated to the weight $\nu$. By
another straightforward calculation one verifies that for any
continuously differentiable function
$g$, setting $\tilde{g}(x)=g(x^2)$, one has
\begin{equation}\label{equation:WishartRewrite}
\langle \nu,(\tilde{g})'V_n\rangle=2Eg'(S+\gamma)U_n(S),
\end{equation}
where  $S$ is a random variable with standard
semicircular law
$\sigma_S$.

\subsubsection{Calculation of $\Var Z_f$ and $E_f$}
\label{subsec-varcalcwish}
For any even polynomial function $f$, 
by the orthogonality relations
noted above, we can write
$$f'(x)=\sum_{n=1}^\infty V_n(x)\,\langle \nu,f'V_n\rangle,$$
with only finitely many nonzero terms in the sum.
Then, by formulas
(\ref{equation:ChebDiag2}) and (\ref{equation:ChebDiag1}), we can in the
present specialization of the band matrix model rewrite
(\ref{equation:CovarianceFormula}) and (\ref{equation:MeanFormula}) in
the form
\begin{equation}\label{equation:CovarianceFormulaWishartRewrite}
\Var Z_{f}=\sum_{r=1}^\infty
(2\lambda_{2r}+\epsilon_{2r})\langle \nu, f'V_r\rangle^2,
 \end{equation}
 \begin{equation}\label{equation:MeanFormulaWishartRewrite}
 E_{f}=\sum_{r=1}^\infty
\frac{1}{2}(\lambda_{2r}+\epsilon_{2r})\langle
\nu, f'V_{2r}\rangle,
 \end{equation}
at least when $f$ is an even polynomial function.
But then these formulas must remain valid for any polynomial $f$ even or
not since $\trace f(X(\NNN_k))$ vanishes
identically for odd $f$.

\subsubsection{Calculation of $\Var Z_{W,g}$}
Given a polynomial function vanishing at the origin consider
the random variables
$$Z_{g,W,k}:=\trace g(W(\NNN_k))-E\trace g(W(\NNN_k))$$
where $g$ is any continuously differentiable function of polynomial
growth. By the bookkeeping
principle (\ref{eq:300304g}), with $\tilde{g}(x)=g(x^2)$, we have 
$Z_{g,W,k}=\frac{1}{2}Z_{\tilde{g},k}$,
hence when $g$ is a polynomial
function, the random variables
$Z_{g,W,k}$ converge in distribution to a mean zero  Gaussian random
variable
$Z_{g,W}$ with variance
\begin{equation}\label{equation:LastFormula}
\Var Z_{g,W}=\sum_{r=1}^\infty (2\lambda_{2r}+\epsilon_{2r})
(Eg'(S+\gamma)U_n(S))^2,
\end{equation}
where, as in formula (\ref{equation:WishartRewrite}), $S$ is a random
variable with law $\sigma_S$. An analogous evaluation of the shift in
mean can also be provided, but to avoid repetitions, we do not state it
here. 

Our final result is

\begin{Theorem}\label{theo-groupwishart}
We work in the setting and under the hypotheses of
Theorems~\ref{Theorem:Preliminary} and \ref{Theorem:Main},
and
in the specialization of the band matrix model
discussed in \S\ref{subsec-secondspec}. If the random variables
$\xi_{\{\alpha,\beta\}}$ satisfy a
\Poincare  inequality with the same constant $c$,
then for any continuously differentiable function $g$ 
with polynomial growth, the random variables $Z_{g,W,k}$ converge in
distribution to a mean zero Gaussian random variable $Z_{g,W}$, with
variance once again given by 
(\ref{equation:LastFormula})\,.
\end{Theorem}
\proof Use
Theorem
\ref{cor-cltpoincare}, for the square-root generalized  Wishart matrices
$X(\NNN_k)$,
taking $\nu$ to be as defined in \S\ref{subsubsection:NuDef},
$\{q_n(x)\}_{n=1}^\infty$ to be the Gram-Schmidt
orthogonalization in $L^2(\nu)$ of the family
 $\{x^{n-1}\}_{n=1}^\infty$ of powers of $x$,
and $f(x)=g(x^2)$.
\qed

\section{Concluding remark}
We have chosen to concentrate in this paper on CLT's for
symmetric matrices. 
Similar techniques work also with Hermitian  matrices,
the main difference being that with
$\xi_{\alpha,\beta}=\xi^1_{\alpha,\beta}+
i\xi_{\alpha,\beta}^2=\xi^*_{\beta,\alpha}$, when $\alpha\neq \beta$,
and $\xi^1_{\alpha,\beta}$ and
$\xi^2_{\alpha,\beta}$ independent, identically
distributed, zero mean real-valued random variables,
it holds that $E[\xi_{\alpha,\beta}]^2=0$, and hence
in the combinatorial
evaluation of the contribution of various terms in expansions
similar to (\ref{equation:WordExpansion2}), the contribution
of words in which an edge is traversed twice 
{\it in the same direction} vanishes. (In particular,
when computing variances for linear statistics
of polynomial type, some of the bracelet contributions
vanish.)
None of the
modifications needed to handle the Hermitian case
are difficult. However, there are sufficiently many such modifications
needed so that 
to give a careful accounting of them would add a nontrivial number of
pages to an already long paper. So we think it best to omit further
discussion.
\\

\noindent
{\it Acknowledgments} We owe the idea to
look at spanning forests when proving 
Lemma \ref{Lemma:Marriage}
to Victor Reiner. We also thank Sergey 
Bobkov for a useful discussion concerning
\Poincare inequalities.

\end{document}